\documentclass[11pt, reqno]{amsart}

\usepackage{amsmath}
\usepackage{amscd}
\usepackage{mathrsfs}
\usepackage{amsfonts}
\usepackage{amssymb}
\usepackage{enumerate}
\usepackage{setspace}
\usepackage{array}

\usepackage{color}

\usepackage[hyperindex, colorlinks=true]{hyperref}

\usepackage[margin=1in]{geometry}

\newcommand{\eqn}{\begin{eqnarray}}
\newcommand{\een}{\end{eqnarray}}

\newcommand{\w}{w dx}

\newcommand{\pat}{\partial_t}
\newcommand{\pax}{\partial_x}
\newcommand{\RR}{\mathbb{R}}

\newtheorem{theorem}{Theorem}[section]

\newtheorem{lemma}{Lemma}[section]

\theoremstyle{definition}
\newtheorem{definition}{Definition}[section]

\newtheorem{remark}{Remark}

\numberwithin{equation}{section}

\begin{document}

\title[Weak solutions to dissipative equations with nonlocal velocity]{Global existence of weak solutions to  dissipative transport equations with nonlocal velocity}

\author{Hantaek Bae}
\address{Department of Mathematical Sciences, Ulsan National Institute of Science and Technology (UNIST), Korea}
\email{hantaek@unist.ac.kr}

\author{Rafael Granero-Belinch\'{o}n}
\address{Univ Lyon, Universit\'e Claude Bernard Lyon 1, CNRS UMR 5208, Institut Camille Jordan, 43 blvd. du 11 novembre 1918, F-69622 Villeurbanne cedex, France.}
\email{granero@math.univ-lyon1.fr}

\author{Omar Lazar}
\address{Instituto de Ciencias Matem\'aticas (ICMAT),  Consejo Superior de Investigaciones Cient\'ificas, Spain}
\email{omar.lazar@icmat.es}

%%%%%%%%%%%General information

\date{\today}



\begin{abstract}
We consider 1D dissipative transport equations with nonlocal velocity field:
\[
\theta_t+u\theta_x+\delta u_{x} \theta+\Lambda^{\gamma}\theta=0, \quad u=\mathcal{N}(\theta),
\]
where $\mathcal{N}$ is a nonlocal operator given by a Fourier multiplier. Especially we consider two types of nonlocal operators:
\begin{enumerate}
\item $\mathcal{N}=\mathcal{H}$, the Hilbert transform,
\item $\mathcal{N}=(1-\partial_{xx} )^{-\alpha}$.
\end{enumerate}
In this paper, we show several global existence of weak solutions depending on the range of $\gamma$, $\delta$ and $\alpha$. When $0<\gamma<1$, we take initial data having finite energy, while we take initial data in weighted function spaces (in the real variables or in the Fourier variables), which have infinite energy, when $\gamma \in (0,2)$.
\end{abstract}

\maketitle

%%%%%%%%%%%%%%%%%%%%
\section{Introduction} \label{sec:1}
%%%%%%%%%%%%%%%%%%%

In this paper, we consider transport equations with nonlocal velocity. Here, the non-locality means that the velocity field is defined through a nonlocal operator that is represented in terms of  a Fourier multiplier. For example, in the two dimensional Euler equation in vorticity form, 
\[
\omega_{t}+u\cdot \nabla \omega=0,
\]
the velocity is recovered from the vorticity $\omega$ through
\[
u=\nabla^\perp(-\Delta)^{-1}\omega \quad \text{or equivalently} \quad  \widehat{u}(\xi)=\frac{i\xi^{\perp}}{|\xi|^{2}}\widehat{\omega}(\xi).
\]
Other nonlocal and quadratically nonlinear equations appear in many applications. Prototypical examples are the surface quasi-geostrophic equation, the incompressible porous medium equation, Stokes equations, magneto-geostrophic equation in multi-dimensions. For more details on nonlocal operators in these equations, see \cite{Bae}.

We here study 1D models of physically important equations. The 1D reduction idea were initiated by Constatin-Lax-Majda \cite{Constantin}: they proposed the following 1D model 
\[
\theta_{t}=\theta\mathcal{H}\theta
\]
for the 3D Euler equation in the vorticity form and proved that $\mathcal{H}\theta$ blows up in finite time under certain conditions. Motivated by this work, other similar models were proposed and analyzed in the literature \cite{Bae, Baker, Carrillo, CC, CC2, Chae, Cordoba 2, De Gregorio, Kiselev, LiRodrigo2, LiRodrigo, Morlet}. In this paper, we consider the following 1D equation:
\eqn \label{model equation}
\theta_t+u\theta_x+\delta u_{x} \theta+\nu\Lambda^{\gamma}\theta=0, \quad u=\mathcal{N}(\theta).
\een
Depending on a nonlocal operator $\mathcal{N}$, (\ref{model equation}) has structural similarity of several important fluid equations as described below. The goal of this paper is to show the existence of weak solutions with rough initial data. To this end, we will choose functionals carefully to extract more information from the structure of the nonlinearity to construct weak solutions.

%%%%%%%%%%%%%%%%%%%%%%%%%%%
\subsection{The case $\mathcal{N}=\mathcal{H}$}
%%%%%%%%%%%%%%%%%%%%%%%%%%%
We first take the case $\mathcal{N}=\mathcal{H}$, the Hilbert transform. Then, (\ref{model equation}) becomes     
\eqn \label{transport equation 1}
\theta_t+\left(\mathcal{H}\theta\right)\theta_x + \delta\theta \Lambda\theta+\nu\Lambda^{\gamma}\theta=0, 
\een
where the range of $\gamma$ and $\delta$ will be specified below. We note that (\ref{transport equation 1}) is considered as an 1D model of the dissipative surface quasi-geostrophic equation. The surface quasi-geostrophic equation describes the dynamics of the mixture of cold and hot air and the fronts between them in 2 dimensions \cite{Constantin2, Pedlosky}. The equation is of the form
\begin{equation}\label{SQG}
\theta_{t}+u\cdot\nabla\theta +\nu\Lambda^{\gamma}\theta=0, \quad u=\left(-\mathcal{R}_{2}\theta, \mathcal{R}_{1}\theta\right),
\end{equation} 
where the scalar function $\theta$ is the potential temperature and $\mathcal{R}_{j}$ is the Riesz transform
\[
\mathcal{R}_{j}f(x)=\frac{1}{2\pi} \text{p.v.} \int_{\mathbb{R}^{2}} \frac{(x_{j}-y_{j})f(y)}{|x-y|^{3}}dy, \quad j=1,2. 
\]
As Constatin-Lax-Majda did for the Euler equation, the equation (\ref{transport equation 1}) is derived by replacing the Riesz transforms with the Hilbert transform. The case $\delta=0$ and $\delta=1$ correspond to (\ref{SQG}) in non-divergence and divergence
form, respectively. We take a parameter $\delta\in [0,1]$ to cover more general nonlinear terms in (\ref{transport equation 1}). We note that there are several singularity formation results when $\nu=0$: $0<\delta<\frac{1}{3}$ and $\delta=1$ \cite{Morlet}, $0<\delta\leq1$  \cite{Chae}, and $\delta=0$ \cite{Cordoba 2, LiRodrigo2, SV}. By contrast, we look for weak solutions of (\ref{transport equation 1}) globally in time (see e.g. \cite{Dong}). From now on, we set $\nu=1$ for notational simplicity.

We now consider (\ref{transport equation 1}). In the direction of seeking a weak solution, we assume that $\theta_{0}$ satisfies the conditions
\eqn \label{initial data 1}
\theta_{0}(x)> 0, \quad  \theta_{0}\in L^{1}\cap H^{\frac{1}{2}}.
\een
Since (\ref{transport equation 1}) satisfies the minimum principle (see Section 2) when $\delta\ge 0$, $\theta(t,x)\ge 0$ for all time. This sign condition combined with the  structure of the nonlinearity enables us to use the following function space
\[
 \mathcal{A}_T= L^{\infty}\left(0,T; L^{p}\cap H^{\frac{1}{2}}\right)\cap L^{2}\left(0,T; H^{\frac{\gamma+1}{2}}\right) \quad \text{for all $p\in (1,\infty)$}.
\]

\textcolor{black}{We note that we add the sign condition of $\theta_{0}$ for the better understanding of the structure of the nonlinearity in the sense that the sign condition and the dissipative term (of any order) deplete the nonlinear effect. If we remove the sign condition of $\theta_{0}$, we would have to control higher other norms, but then we will obtain strong solutions locally in time  which are not the notion of solutions we are going to obtain in this paper. Note that the positivity condition is something that is also assumed in other works (see e.g. \cite{Cordoba 2}, \cite{Dong}, \cite{Kiselev}). }

\begin{definition}
We say $\theta$ is a weak solution of (\ref{transport equation 1}) if $\theta\in \mathcal{A}_T$ and (\ref{transport equation 1}) holds in the following sense:  for any test function $\psi\in \mathcal{C}^{\infty}_{c}\left([0,T)\times\mathbb{R}\right)$, 
\[
\int^{T}_{0}\int_{\mathbb{R}} \left[\theta \psi_{t} + \left(\mathcal{H}\theta\right)\theta\psi_{x}+(1-\delta)\Lambda\theta\theta \psi -\theta \Lambda^{\gamma}\psi\right]  dxdt =\int_{\mathbb{R}}\theta_{0}(x)\psi(0,x)dx
\]
holds for any $0<T<\infty$.
\end{definition}

\begin{theorem}\label{GW 1}
Let $\gamma\in (0,1)$ and $\delta\ge \frac{1}{2}$. Then, for any  $\theta_{0}$ satisfying (\ref{initial data 1}), there exists a weak solution of (\ref{transport equation 1}) in $\mathcal{A}_{T}$ for all $T>0$. Moreover, a weak solution is unique when $\gamma=1$.
\end{theorem}

For $\gamma \in (0,2)$, we consider infinite energy solutions of (\ref{transport equation 1}). More precisely, we take a family of weights $w_{\beta}=\left(1+|x|^{2}\right)^{-\frac{\beta}{2}}$, $0<\beta<\gamma$ and we shall prove various existence theorems. For the critical case, we take initial data in the following weighted Sobolev space  
%\eqn \label{infinite initial data}
%\theta_{0}(x)\ge 0, \quad  \theta_{0}\in H^{\frac{1}{2}}\left(w_{\beta}dx\right)\cap L^{\infty}
%\een
\eqn \label{infinite initial data}
 \theta_{0}\in H^{\frac{1}{2}}\left(w_{\beta}dx\right)\cap L^{\infty}.
\een
These weighted spaces are defined in Section 2. We note that $\theta_{0}$ can decay (slowly) at infinity. For example, as long as $\beta-2\eta>1$, $|\theta_{0}(x)|\simeq |x|^{-\eta}$, $\eta>0$ is allowed to stay in $L^{2}(w_{\beta}dx)$
\[
\int_{|x|\ge 1}\frac{|x|^{2\eta}}{(1+|x|^{2})^{\frac{\beta}{2}}}dx<\infty.
\]
But, we can still use the energy method to obtain a weak solution of (\ref{transport equation 1}). Let
\[
\mathcal{B}_T= L^{\infty}\left(0,T; H^{\frac{1}{2}}(w_{\beta}dx)\right)\cap L^{2}\left(0,T; H^{1}(w_{\beta}dx)\right).
\]

\begin{theorem}\label{GW infinite energy}
Assume $\gamma=1$, then, for any  $\theta_{0}$ satisfying (\ref{infinite initial data}) with $\|\theta_{0}\|_{L^{\infty}}$ being sufficiently small, there exists a unique weak solution of (\ref{transport equation 1}) in $\mathcal{B}_{T}$ for all $T>0$. 
\end{theorem}
{\color{black}\begin{remark}It is worth mentioning that this theorem is also true in the unweighted setting, this is an important point because one would need to use the existence of solutions in the unweighted setting to prove Th. \ref{GW infinite energy}. Indeed, one has nice {\it{a priori}} estimates in the unweighted setting, it suffices to observe that the evolution of the $\dot H^{1/2}$ (semi-)norm is obtained via the classical Hardy-BMO duality along with formula \ref{Hilbert identity}. Indeed, one writes 
\begin{eqnarray*}
\frac{1}{2} \frac{d}{dt} \int \Vert\Lambda^{1/2} \theta \Vert^{2}_{L^{2}} \ dx + \int \vert  \Lambda \theta \vert^{2} \ dx &=& -\int \Lambda \theta \mathcal{H} \theta \theta_{x} \ dx  - \int \theta \vert \Lambda \theta \vert^{2} \ dx \\
&\leq& \Vert \theta_{x} \mathcal{H} \theta_{x} \Vert_{\mathcal{H}^{1}} \Vert \mathcal{H}\theta \Vert_{BMO} + \Vert \theta_{0} \Vert_{L^{\infty}} \Vert  \Lambda \theta \Vert^{2}_{L^{2}} \\
&\leq&2\Vert \theta_{x} \Vert^{2}_{L^{2}} \Vert \theta_{0} \Vert_{L^{\infty}}=2\Vert \mathcal{H}\theta_{x} \Vert^{2}_{L^{2}} \Vert \theta_{0} \Vert_{L^{\infty}}.
\end{eqnarray*}
\end{remark}
\noindent Hence, if $\Vert \theta_{0} \Vert_{L^{\infty}}<1/2$, one has $\theta \in \mathcal{C}([0,T], \dot H^{1/2}) \cap L^{2}( [0,T], \dot H^{1})$ for all finite $T>0$. For the construction using compactness we refer to \cite{Lazar}.} \\

In the subcritical case i.e. $\gamma \in (1,2)$, we have global existence of weak solutions for any arbitrary initial data in the weighted Sobolev space $H^{1}( \w)$, with $w(x)=w_{\beta}(x)=(1+x^2)^{-\beta/2}, \beta \in (0,1)$.
\begin{theorem} \label{sub}
 Let $\gamma \in (1,2)$, for all $\theta_0 \in H^{1}( \w) \cap L^\infty$  there exists at least one global weak solution to the equation $\mathcal{T}_{\alpha},$ which verifies, for all finite $T>0$
$$
\theta \in \mathcal{C}([0,T], H^{1}( \w)) \cap L^{2}([0,T] , \dot H^{1+\alpha/2}( \w)).
$$
Moreover, for all $T<\infty$, we have
$$\Vert \theta (T) \Vert^{2}_{H^{1}( \w)} \leq \Vert \theta_{0} \Vert^{2}_{H^{1}( \w)} e^{CT} $$ \\
\vspace{0cm}
\noindent The constant $C>0$ depends  only on  $\Vert \theta_{0} \Vert_{L^\infty}$, $\beta$, $k$, $\delta$ and $\nu$. \\

\end{theorem} 
In the supercritical case, one can prove the following local existence theorem for data in the weighted $H^{2} (\w)$ Sobolev space.

\begin{theorem} \label{super}   Assume that $0<\alpha<1$ and $\delta\geq0$, then for all positif initial data  $\theta_{0} \in H^{2}( \w)$ where the weight is given by $w_{\beta}(x)=(1+\vert x\vert^{2})^{-\beta/2}$ with $\beta\in (0,\alpha/2)$, there exists a time $T^{*}(\theta_0)>0$ such that  $(\mathcal{T}_{\alpha})$ admits at least one solution that verifies $$\theta \in \mathcal{C}([0,T], H^{2}( \w)) \cap L^{2}([0,T], \dot H^{2+\frac{\alpha}{2}}( \w))$$ for all $T\leq T^{*}$.
\end{theorem}

We have some restrictions on the sign of initial data and the range of $\delta$ in Theorem \ref{GW 1}, and the smallness condition of $\|\theta_{0}\|_{L^{\infty}}$ in Theorem \ref{GW infinite energy}.   We can remove these conditions by looking for a solution of (\ref{transport equation 1}) in function spaces defined by the Fourier transform. Let 
\[
A^\alpha=\left\{f\in L^1_{loc}:  \|f\|_{A^\alpha}=\int_{\RR} (1+|\xi|^\alpha)|\hat{f}(\xi)|d\xi <\infty\right\}.
\]
We also define 
\[
\mathcal{W}_T= L^{\infty}\left(0,T; W^{1,\infty}\right) \cap W^{1,\infty}(0,T;L^\infty)\cap L^1(0,T;W^{2,\infty}).
\]

\begin{theorem}\label{GW Wiener}
Let $\gamma=1$ and $\delta\in\RR$. Then, for any $\theta_{0}\in A^1$ with 
\eqn \label{initial data in Wiener}
\|\theta_{0}\|_{A^0}<\frac{\sqrt{\pi}}{\sqrt{2}(1+|\delta|)},
\een
there exists a unique weak solution of (\ref{transport equation 1}) verifying the following inequality for all $T>0$
\[
\theta\in\mathcal{W}_T, \quad \sup_{t\in[0,T]}\|\theta(t)\|_{A^1}+\left(1-\frac{\sqrt{2} (1+|\delta|)\|\theta_{0}\|_{A^{0}}}{\sqrt{\pi}}\right)\int^{T}_{0}\|\theta_{x}(t)\|_{A^1} dt \leq \|\theta_0\|_{A^1}.
\]
\end{theorem}

We note that $\theta_{0}\in A^{1}$ can have infinite energy. For example, we take $\widehat{\theta_{0}}(\xi)=\frac{e^{-|\xi|}}{\sqrt{|\xi|}}$ for $\xi\ne 0$. Then, $\theta_{0}\in A^{1}$ but $\theta_{0}\notin L^{2}$.

\textcolor{black}{We observe that the proof of Theorem \ref{GW Wiener} is due to the perfect balance (in the critical case $\gamma=1$) between the derivatives in the nonlinearity and the diffusive linear operator. This is due to the fact that, in Wiener spaces, the parabolic gain of regularity for $\Lambda$ is $L^1_tA_x^1$ (\emph{i.e.} a full derivative in the Wiener space). This is in contrast with the case of $L^2$, where the parabolic gain of regularity for $\Lambda$ is $L^2_tH_x^{1/2}$ (\emph{i.e.} just half derivative in $L^2$). In particular, the proof is based in an inequality of the type
$$
\frac{d}{dt}E(t)+E(t)D(t)\leq -D(t),
$$
where $E$ and $D$ are the appropriate energy and dissipation. From such an inequality the decay of the energy for small enough initial energy can be easily obtained. In the case where $\gamma>1$, that balance is broken and the previous inequality has to be replaced by
$$
\frac{d}{dt}E(t)+E^{2-\alpha}(t)D^{\alpha}(t)\leq -D(t),
$$
with $0<\alpha<1$. In our setting, the lack of Poincar\'e inequality that relates $E$ and $D$, makes this new inequality less suited for our approach.
}

%%%%%%%%%%%%%%%%%%%%%%%%%%%%%%%%%%%%%%%%%
\subsection{The case $\mathcal{N}=(1-\partial_{xx} )^{-\alpha}$ and $\delta=0$}
%%%%%%%%%%%%%%%%%%%%%%%%%%%%%%%%%%%%%%%%%
In this case, (\ref{model equation}) is changed to the equation
\eqn \label{transport equation 2}
\theta_{t}+u\theta_{x}+ \delta\theta \Lambda\theta+\Lambda^{\gamma}\theta=0, \quad u=(1-\partial_{xx} )^{-\alpha} \theta.
\een
This equation is closely related to a generalized Proudman-Johnson equation \cite{Okamoto, Proudman, Wunsch}:
\[
f_{txx} + f f_{xxx} + \delta f_{x}f_{xx} = \nu f_{xxxx}
\]
which is derived from the 2D incompressible Navier-Stokes equations via the separation of space variables when $\delta=1$. By taking $w=f_{xx}$, 
\[
w_t +fw_x +\delta f_x w=\nu w_{xx}, \quad f=(\partial_{xx})^{-1}w.
\]
The inviscid case with $\delta=2$ is equivalent to the Hunter-Saxton equation arising in the study of nematic liquid crystals \cite{Hunter}. The equation (\ref{transport equation 2}) is also considered as a model equation of the Lagrangian averaged Navier-Stokes  equations \cite{Marsden} which are given by
\[
\partial_{t}\left(1-\sigma^{2}\Delta\right)u+u\cdot \nabla \left(1-\sigma^{2}\Delta\right)u +(\nabla u)^{T} \cdot \left(1-\sigma^{2}\Delta\right)u=-\nabla p +\nu \Delta \left(1-\sigma^{2}\Delta\right)u, \quad \nabla \cdot u=0.
\]

\textcolor{black}{
But, we here consider (\ref{transport equation 2}) with $\delta=0$ to see how the regularizing effect in $u$ overcomes difficulties from the quadratic term $u\theta_{x}$. Along this direction, we look closely to see if there is a solution when $\alpha$ and $\gamma$ meet certain conditions.}

We first deal with (\ref{transport equation 2}) with initial data in $L^{2}\cap L^{\infty}$. Let 
\[
\mathcal{C}_T= L^{\infty}\left(0,T; L^{p}\right)\cap L^{2}\left(0,T; H^{\frac{\gamma}{2}}\right) \quad \text{for all $p\in [2,\infty]$}.
\]

\begin{definition}
We say $\theta$ is a weak solution of (\ref{transport equation 2}) if $\theta\in \mathcal{C}_T$ and (\ref{transport equation 2}) holds in the following sense:  for any test function $\psi\in \mathcal{C}^{\infty}_{c}\left([0,T)\times\mathbb{R}\right)$, 
\[
\int^{T}_{0}\int_{\mathbb{R}} \left[\theta \psi_{t} +u_{x}\theta \psi+ u\theta \psi_{x}-\theta \Lambda^{\gamma}\psi\right]  dxdt =\int_{\mathbb{R}}\theta_{0}(x)\psi(0,x)dx
\]
holds for any $0<T<\infty$.
\end{definition}

\begin{theorem} \label{GW 2}
Let $\gamma\in (0,2)$ and $\alpha= \frac{1}{2}-\frac{\gamma}{4}$.  Then, for any $\theta_0\in L^{2}\cap L^{\infty}$, there exists a weak solution of (\ref{transport equation 2}) in $\mathcal{C}_{T}$ for all $T>0$. \textcolor{black}{Moreover, such a weak solution is unique if $\gamma>1$ and $\alpha>\frac{1}{4}$.}
\end{theorem}

\textcolor{black}{We note that $\theta_0\in L^{2}\cap L^{\infty}$ is enough to construct a weak solution in Theorem \ref{GW 2}, but we need to strengthen $\gamma$ and $\alpha$ to get the uniqueness of weak solutions.  }

When $\gamma=1$, we consider weights $w_{\beta}(x)=\left(1+|x|^{2}\right)^{-\frac{\beta}{2}}$ with $0<\beta<1$, and take initial data in $H^{1}(w_{\beta}dx)\cap L^{\infty}$. Let $\alpha=\frac{1}{4}$ and 
\[ 
\mathcal{D}_T= L^{\infty}\left(0,T; H^{1}(w_{\beta}dx)\right)\cap L^{2}\left(0,T; H^{\frac{3}{2}}(w_{\beta}dx)\right).
\]

\begin{theorem}\label{GW infinite energy 2}
Let $\gamma=1$ and $\alpha=\frac{1}{4}$. Then, for any  $\theta_{0}\in H^{1}(w_{\beta}dx)\cap L^{\infty}$, there exists a unique global weak solution of (\ref{transport equation 2}) in $\mathcal{D}_{T}$ for all $T>0$. 
\end{theorem}

Compared to Theorem \ref{GW infinite energy}, we do not assume that $\|\theta_{0}\|_{L^{\infty}}$ is small to prove Theorem \ref{GW infinite energy 2}.

\vspace{1ex}

In Theorem \ref{GW 2} and Theorem \ref{GW infinite energy 2}, we have restrictions on the range of $\alpha$. Again, we can remove these conditions by looking for a solution of (\ref{transport equation 2}) in function spaces defined by the Fourier variables.

\begin{theorem}\label{GW Wiener2}
Let $\gamma=1$ and $\alpha\geq0$. Then, for any $\theta_{0}\in A^1$ satisfying 
\eqn  \label{initial data in Wiener 2}
\|\theta_{0}\|_{A^0}<\frac{\sqrt{\pi}}{\sqrt{2}},
\een
there exists a unique weak solution of (\ref{transport equation 2}) verifying the following inequality for all $T>0$
\[
\theta\in\mathcal{W}_T, \quad \sup_{t\in[0,T]}\|\theta(t)\|_{A^1}+\left(1-\frac{\sqrt{2} \|\theta_{0}\|_{A^{0}}}{\sqrt{\pi}}\right)\int^{T}_{0}\|\theta_{x}(t)\|_{A^1} dt\leq \|\theta_0\|_{A^1}.
\]
\end{theorem}

\begin{remark}
We note that Theorem \ref{GW Wiener2} remains valid with straightforward changes in the spirit of Theorem \ref{GW Wiener} when $\delta\neq0$.
\end{remark}

%%%%%%%%%%%%%
\section{Preliminaries}
%%%%%%%%%%%%%
All constants will be denoted by $C$ that is a generic constant. In a series of inequalities, the value of $C$ can vary with each inequality. For $s\in \mathbb{R}$, $H^{s}$ is a Hilbert space with 
\[
\|f\|^{2}_{H^{s}}=\int_{\mathbb{R}}(1+|\xi|^{2})^{s}\left|\hat{f}(\xi)\right|^{2}d\xi.
\]
%%%%%%%%%%%%%%%%%%%%%%%%%%%%%%%%%%%
\subsection{Hilbert transform and fractional Laplacian}
%%%%%%%%%%%%%%%%%%%%%%%%%%%%%%%%%%
The Hilbert transform is defined as 
\[
\mathcal{H}f(x)=\text{p.v.} \int_{\mathbb{R}} \frac{f(y)}{x-y}dy.
\]
The differential operator $\Lambda^{\gamma}=(\sqrt{-\Delta})^{\gamma}$ is defined by the action of the following kernels \cite{Cordoba}:
\eqn \label{lambda gamma}
\Lambda^{\gamma} f(x)=c_{\gamma}\text{p.v.} \int_{\mathbb{R}} \frac{f(x)-f(y)}{|x-y|^{1+\gamma}}dy,
\een
where $c_{\gamma}>0$ is a normalized constant. When $\gamma=1$,   
\[
\Lambda f(x)=\mathcal{H}f_{x}(x).
\]

Moreover, we have the following identity:
\eqn \label{Hilbert identity}
\mathcal{H} \left(\theta_{x}\left(\mathcal{H}\theta_{x}\right)\right) =\frac{1}{2}\left[\left(\Lambda \theta\right)^{2}-\left(\theta_{x}\right)^{2}\right].
\een

We also recall the following pointwise property of $\Lambda^{\alpha}$.

\begin{lemma}\label{pointwise lemma}\cite{Cordoba}
Let $0\leq \alpha \leq 2$ and $f\in \mathcal{S}$. Then, 
\begin{equation*}
\begin{split}
&f(x)\Lambda^{\alpha} f(x) \ge \frac{1}{2}\Lambda^{\alpha} \left(f^{2}(x)\right), \\
& f^{2}(x)\Lambda f(x) \ge \frac{1}{3}\Lambda \left(f^{3}(x)\right) \ \text{when $f\ge 0$} 
\end{split}
\end{equation*}
\end{lemma}

%%%%%%%%%%%%%%%%%%%%%%%%%%%%
\subsection{Minimum and Maximum Principles}
%%%%%%%%%%%%%%%%%%%%%%%%%%%%
In Theorem \ref{GW 1}, we assume $\theta_{0}>0$. To obtain global-in-time solutions, we need $\theta(t,x)\ge 0$ for all time. We first assume that $\theta(t,x)\in C^1\left([0,T]\times\mathbb{R}\right)$ and $x_t$ be a point such that $\displaystyle m(t)=\theta(t,x_t)$. If $m(t)>0$ for all time, nothing is left to prove. So, we check a point $(t, x_{t})$ where $m(t)=0$. Since $m(t)$ is a continuous Lipschitz function, it is differentiable at almost every $t$ by Rademacher's theorem. From the definition of $\Lambda^{\gamma}$, 
\begin{equation*}
\begin{split}
\frac{d}{dt}m(t)&=-\delta \theta(t,x_{t})\text{p.v.} \int_{\mathbb{R}} \frac{\theta(t,x_{t})-\theta(t,y)}{|x_{t}-y|^{1+\gamma}}dy - \text{p.v.} \int_{\mathbb{R}} \frac{\theta(t,x_{t})-\theta(t,y)}{|x_{t}-y|^{1+\gamma}}dy\\
& \ge \left[-\delta \text{p.v.} \int_{\mathbb{R}} \frac{\theta(t,x_{t})-\theta(t,y)}{|x_{t}-y|^{1+\gamma}}dy\right]m(t).
\end{split}
\end{equation*}
Since the quantity in the bracket is nonnegative when $\delta\ge 0$, we have that $m(t)$ is non-decreasing in time if $\theta_{0}>0$ and thus $\theta(t,x)\ge 0$ for all time. Similarly, maximum values of $\theta(t,x)$ are non-increasing in time when $\theta_{0}>0$ with $\theta_{0}\in L^{\infty}$. For general initial data satisfying (\ref{initial data 1}) and (\ref{infinite initial data}), we can use regularization method. For such a regularized problem with smooth solution $\theta^\epsilon$, the same argument works. Then, we construct $\theta$ as the limit of $\theta^\epsilon$. As $\theta$ will be also the pointwise limit of $\theta^\epsilon$ almost everywhere, we conclude that $\theta(t,x)\ge 0$.

Since (\ref{transport equation 2}) is purely a dissipative transport equation, we immediately have 
\[
\|\theta(t)\|_{L^{\infty}}\leq \|\theta_{0}\|_{L^{\infty}}.
\]

%%%%%%%%%%%%%%%%%%%%%%%%%%%%%%%%%%%
\subsection{The Wiener spaces $A^\alpha$}
%%%%%%%%%%%%%%%%%%%%%%%%%%%%%%%%%%
The Wiener space is defined as
\[
A^0=\left\{f\in L^1_{loc}: \hat{f}(\xi)\in L^1\right\},
\]
where $\widehat{f}$ denotes the Fourier transform of $f$ 
\[
\widehat{f}(\xi)=\frac{1}{\sqrt{2\pi}}\int_{\RR}f(x)e^{-ix\cdot\xi}dx.
\]
$A^0$ is a Banach space endowed with the norm
\[
\|f\|_{A^0}=\|\widehat{f}\|_{L^1}.
\]
Furthermore, using Fubini's Theorem, $A^0$ is a Banach algebra, \emph{i.e.} 
\[
\left\|fg\right\|_{A^0}\leq \|f\|_{A^0}\|g\|_{A^0}.
\]
Once we have defined $A^0$, we can define the full scale of homogeneous, $\dot{A}^\alpha$, and inhomogeneous, $A^\alpha$, Wiener spaces as
\begin{equation}\label{eq:4b}
\begin{split}
&\dot{A}^\alpha=\left\{f\in L^1_{loc}:  \|f\|_{\dot{A}^\alpha}=\int_{\RR} |\xi|^\alpha|\hat{f}(\xi)|d\xi<\infty\right\},\\
& A^\alpha=\left\{f\in L^1_{loc}: \|f\|_{A^\alpha}=\int_{\RR} (1+|\xi|^\alpha)|\hat{f}(\xi)|d\xi
\right\}.
\end{split}
\end{equation}
For these spaces, the following inequalities hold
\begin{align}\label{eq:5a}
\|f\|_{C(\RR)}&\leq \|f\|_{A^0(\RR^d)} \quad \forall\, f\in A^0(\RR)\\
\label{eq:5d}
\|f\|_{\dot{A}^\alpha(\RR)}&\leq \|f\|_{A^0(\RR^d)}^{1-\theta}\|f\|_{\dot{A}^{\frac{\alpha}{\theta}}(\RR^d)}^{\theta}  \quad \forall\, 0<\theta<1, \ \alpha\geq0,\ f\in A^0(\RR)\cap A^{\frac{\alpha}{\theta}}(\RR).
\end{align}

As a consequence of \eqref{eq:5a}, we obtain that if $u\in A^0$ has infinite energy then
$$
\limsup_{|x|\rightarrow\infty} |u(x)|+\liminf_{|x|\rightarrow\infty} |u(x)|<\infty.
$$

%%%%%%%%%%%%%%%%%%%%
\subsection{Commutator estimate}
%%%%%%%%%%%%%%%%%%%%
In the proof of Theorem \ref{GW 1}, we need to estimate a commutator term involving $\Lambda^{\frac{1}{2}}$. To do this, we first recall Hardy-Littlewood-Sobolev inequality in 1D. Let $K_{\alpha}(x)=\frac{1}{|x|^{\lambda}}$ and $T_{\lambda}f=K_{\lambda}\ast f$. Then, 
\[
\left\|T_{\lambda}f\right\|_{L^{q}}\leq C \|f\|_{L^{p}}, \quad \frac{1}{q}+1=\frac{1}{p}+\lambda.
\]

\begin{lemma} \label{commutator lemma}
For $f\in L^{\frac{3}{2}}$, $g\in L^{\frac{3}{2}}$  and $\psi \in W^{1,\infty}$, 
\[
\left\| \left[\Lambda^{\frac{1}{2}},\psi \right]f- \left[\Lambda^{\frac{1}{2}},\psi \right]g\right\|_{L^{6}}\leq C\|\psi\|_{W^{1,\infty}} \left\|f-g\right\|_{L^{\frac{3}{2}}}.
\]
\end{lemma}

\begin{proof}
By the definition of $\Lambda^{\frac{1}{2}}$, we have
\[
\left(\left[\Lambda^{\frac{1}{2}},\psi \right]f- \left[\Lambda^{\frac{1}{2}},\psi \right]g\right)(x)=c_{1}\text{p.v.}\int \frac{(\psi(y)-\psi(x))(f(y)-g(y))}{|x-y|^{\frac{3}{2}}}dy
\]
and thus
\begin{equation}
\begin{split}
\left|\left[\Lambda^{\frac{1}{2}},\psi \right]f- \left[\Lambda^{\frac{1}{2}},\psi \right]g\right|(x) \leq C\|\nabla \psi\|_{L^{\infty}} \int \frac{|f(y)-g(y)|}{|x-y|^{\frac{1}{2}}}dy.
\end{split}
\end{equation}
Using Hardy-Littlewood-Sobolev inequality, we obtain that
\begin{equation}
\begin{split}
\left\|\left[\Lambda^{\frac{1}{2}},\psi \right]f- \left[\Lambda^{\frac{1}{2}},\psi \right]g\right\|_{L^{6}} \leq C\left\|\nabla \psi\right\|_{L^{\infty}}  \left\|f-g\right\|_{L^{\frac{3}{2}}}
\end{split}
\end{equation}
which completes the proof.
\end{proof}

%%%%%%%%%%%%%%%%%%%%%
\subsection{Muckenhoupt weights}
%%%%%%%%%%%%%%%%%%%%
We briefly introduce weighted spaces. A weight  $w$ is a positive and locally integrable function. A measurable function $\theta$ on $\mathbb{R}$ belongs to the weighted Lebesgue spaces $L^{p}(wdx)$ with $1\leq p<\infty$ if and only if
\[
\|\theta\|^{p}_{L^{p}(wdx)}=\int_{\mathbb{R}} |\theta(x)|^{p}w(x)dx<\infty.
\]
An important class of weights is the Muckenhoupt class $\mathcal{A}_{p}$ for $1<p<\infty$ \cite{Coifman, Muckenhoupt}. Let $1<p<\infty$, we say that $w \in \mathcal{A}_{p}$ if and only if there exists a constant $C_{p,w}>0$ such that
\[
\sup_{r>0,x_{0} \in \mathbb{R}} \left( \frac{1}{2r} \int_{[x_{0}-r,x_{0}+r]} w dx \right)\left( \frac{1}{2r} \int_{[x_{0}-r,x_{0}+r]} w^{\frac{1}{1-p}} dx \right)^{p-1} \leq C_{p,w}.
\]

 This class satisfies the following properties.
\begin{enumerate}[]
\item (1) Calder\'on-Zygmund type operators are bound on $L^{p}(wdx)$ when $w\in \mathcal{A}_{p}$ and $1<p<\infty$ \cite{Stein}.
\item (2) Let $w\in \mathcal{A}_{p}$. We define weighted Sobolev spaces as follows
\begin{equation} \label{weighted Sobolev space}
\begin{split}
& f\in H^{1}(wdx) \iff f\in L^{2}(wdx)  \ \text{and} \ f_{x}\in L^{2}(wdx),\\
& f\in H^{1}(wdx) \iff (1-\partial_{xx})^{\frac{1}{2}}f\in L^{2}(wdx) \iff f \in L^{2}(wdx)  \ \text{and} \ \Lambda f\in L^{2}(wdx),\\ 
& f\in H^{\frac{1}{2}}(wdx) \iff (1-\partial_{xx})^{\frac{1}{4}}f\in L^{2}(wdx) \iff f \in L^{2}(wdx)  \ \text{and} \ \Lambda^{\frac{1}{2}}f\in L^{2}(wdx).
\end{split}
\end{equation}
\item (3) Gagliardo-Nirenberg type inequalities (see e.g \cite{Maz})
\begin{equation}\label{GN weight}
\begin{split}
& \left\|\Lambda^{\frac{1}{2}}f\right\|_{L^{2}(wdx)}\leq C \left\|f\right\|^{\frac{1}{2}}_{L^{2}(wdx)}  \left\|\Lambda f\right\|^{\frac{1}{2}}_{L^{2}(wdx)},\\
& \|\theta\|_{L^{4}(wdx)}\leq C \|\theta\|^{\frac{1}{2}}_{L^{2}(wdx)} \left\|\Lambda^{\frac{1}{2}}\theta\right\|^{\frac{1}{2}}_{L^{2}(wdx)}.
\end{split}
\end{equation}
\end{enumerate}
This latter inequality can be proved for instance by using the weighted Sobolev embedding $H^{\frac{1}{4}}(wdx) \hookrightarrow L^{4}(wdx)$, and then by weighted interpolation one recover the second inequality in (\ref{GN weight}).  

In this paper, we take weights $w_{\beta}=(1+|x|^{2})^{-\frac{\beta}{2}}$, $0<\beta<1$, which belongs to the $\mathcal{A}_{p}$ class of Muckenhoupt for all $1<p<\infty$. These weights also satisfy the following properties. For the proofs, see \cite{Lazar} (for the first two points), and \cite{Lazar2} (for the last two points).

\begin{lemma}\label{lem24}
Let $w_{\beta}(x)=(1+|x|^{2})^{-\frac{\beta}{2}}$. 
\begin{itemize}
\item Let $p\geq2$ be such that $\frac{3}{2}-\beta(1-\frac{1}{p}) >1$, then the commutator $\frac{1}{w_{\beta}} [\Lambda^{1/2},w_{\beta}]$ is bounded from $L^{p}(\w)$ to $L^{p}(\w)$. 
 \item Let $2\leq p < \infty$,  then the commutator $\frac{1}{\sqrt{w}} [\Lambda,\sqrt{w}]$ is bounded from $L^{p}(\w)$ to $L^{p}(\w )$.
\item  Let $s \in (0,1)$, then for all $\beta \in(0,s)$,  the commutator $\frac{1}{{w}} [\Lambda^{s/2},{w}]$ is bounded from $L^{2}(\w )$ to $L^{2}(\w)$ 
\item  Let $s \in (1,2)$, then for all $\beta \in (0,1)$,  the commutator $\frac{1}{{w}} [\Lambda^{s/2},{w}]$ is  continuous from $L^{2}(\w)$ to $L^{2}(\w).$
\item The commutator $\frac{1}{w} [\Lambda^{1/2},w]$ is bounded from $L^{p}(\w)$ to $L^{p}(\w)$
\end{itemize}
\end{lemma}
\begin{remark} Actually, the last point of \ref{lem24} is not proved in  \cite{Lazar} or  \cite{Lazar2}. However, a slight modification of the proof given in \cite{Lazar} of the second point of \ref{lem24} gives the last point of \ref{lem24}. Indeed, the idea of the proof of all of those commutator estimates is to split the domain of integration into 3 regions ($\Delta_{1}, \Delta_{2}, and  \Delta_{3}$ according to the notation in  \cite{Lazar}). Since we are dealing with a quite singular kernel, we need to do a sort of second Taylor expansion in the region  where the kernel is very singular (the region $\Delta_{1}$), the same argument works writting $w$ in stead of $\gamma$. As well, in the region $\Delta_{2}$ one can follow \cite{Lazar}. In the region $\Delta_{3}$, the estimates of the kernel change a little bit, indeed, we have the following estimate
$$
\vert K(x,y) \vert \leq C \frac{w(x)^{\frac{1}{p}-1} + w(y)^{-\frac{1}{p}}}{\vert x-y \vert^{3/2}} \leq   \frac{C'}{\vert x-y \vert^{2-\beta(1-\frac{1}{p})}} +  \frac{C'}{\vert x-y \vert^{2-\frac{\beta}{p}}} \leq \frac{C'}{\vert x-y \vert^{2-\beta \max(1-\frac{1}{p},\frac{1}{p})}} 
$$
where we used that, on $\Delta_3$ we have $1 \leq w(x)^{-1} \leq C \vert x-y \vert^{\beta}$ and $1 \leq w(y)^{-1} \leq C \vert x-y \vert^{\beta}$. To conclude, it suffices to observe that since $0<\beta<1$ and $\max(1-\frac{1}{p},\frac{1}{p}) \in (0,1]$  for all $p\geq2$, we have  $K \in L^{1}$. 
\end{remark}

We shall also need to estimate $ \vert \Lambda^{s} w \vert$ for $s\in (0,2)$, the estimate of such a term will be done using the following lemma (see \cite{Lazar} and \cite{Lazar2} for the proof).
\begin{lemma}
Let $w_{\beta}(x)=(1+|x|^{2})^{-\frac{\beta}{2}}$, with $\beta \in (0,1)$.
\begin{itemize}
\item $\vert \partial_{x} w_{\beta}(x) \vert \leq C(\beta) w_{\beta}(x)$
\item For all $s \in (0,2)$,  there exists a  constant $C=C(s,\beta)>0$ such that the following estimate holds
$$\vert \Lambda^{s} w_{\beta}(x) \vert \leq C  w_{\beta}(x)$$
\end{itemize}
\end{lemma}

\subsection{Littlewood-Paley decomposition} 

Consider a fixed function $\phi_{0} \in \mathcal{D}(\mathbb R)$ that is non-negative and radial and is such that ${\phi}_{0}(\xi) = 1$,  if $\vert \xi \vert \leq 1/2$  and ${\phi}_{0}(\xi)=0$  if $\vert \xi \vert \geq 1$. Then, we define a new function ${\psi}_{0}$ :  ${\psi}_{0}(\xi)=\phi_{0}(\xi/2)-\phi_{0}(\xi)$ (which is supported in a corona). Then, for $j \in\mathbb{Z}$, we define the two distributions $S_{j}f= \mathcal{F}^{-1}(\phi_{0}(2^{-j}\xi ) \hat{f}(\xi))$ and
 $\Delta_{j}f= \mathcal{F}^{-1}(\psi_{0}( 2^{-j}\xi) \hat{f}(\xi))$ and we get the so-called inhomogeneous Littlewood-Paley decompositon of 
$f \in \mathcal{S}'(\mathbb R)$ that is for all $K\in \mathbb Z$  the following inequality holds in $\mathcal{S}'(\mathbb R)$
\begin{equation} \label{lp}
f=S_{K} f + \sum_{j\geq K} \Delta_{j}f.
\end{equation} 
Passing to the limit in equality \ref{lp} as $K\rightarrow -\infty$ in the $\mathcal{S}'(\mathbb R)$ topology one obtains  
\begin{equation} \label{lph}
f= \sum_{j\in \mathbb Z} \Delta_{j}f.
\end{equation} 
The equality \ref{lph} is called homogeneous decomposition of $f$ and is defined modulo polynomials. We are now ready to define the homogeneous weighted Sobolev spaces $\dot H^{s}_{w}$ for $\vert s \vert <1/2$, they are defined as follows
$$
f \in \dot H^{s}_{w}\Longleftrightarrow  f =\sum_{j\in \mathbb Z} \Delta_{j}f,   \ \ {\text {in}} \  \mathcal{S}'(\mathbb R)\ \ {\text{and}} \ \  \sum_{j\in \mathbb Z} 2^{2j} \Vert \Delta_{j}f  \Vert^{2}_{L^{2}_{w}}<\infty.
$$
We shall use the  Bernstein's inequality, that is for all $f \in \mathcal{S}'(\mathbb{R})$ and   $(j, s) \in \mathbb{Z} \times \mathbb{R}$, and for all $1\leq p \leq q \leq \infty$ and all weights $w \in A_{\infty}$, we have 
 $$\Vert \Lambda^{s} \Delta_{j} f \Vert_{L^{p}_{w}} \lesssim 2^{js} \Vert \Delta_j f \Vert_{L^{p}_{w}} \ \ {\text{also}} \ \ 
 \Vert \Delta_j f \Vert_{L^{q}_{w}} \lesssim2^{{j}(\frac{1}{p}-\frac{1}{q})}  \Vert \Delta_j f \Vert_{L^{p}_{w}} \ \text{and}  \  
\Vert \Lambda^{s} S_j f \Vert_{L^{p}_{w}} \lesssim 2^{js} \Vert S_j f \Vert_{L^{p}_{w}} $$
Finally, we recall that for any distributions $f$ and $g$ that are in $\mathcal{S}'(\mathbb R)$, one has the following paraproduct formula 
\begin{equation} \label{ppf}
fg=\sum_{q \in \mathbb Z} S_{q+1} f \Delta_{q} g + \sum_{j \in \mathbb Z} \Delta_{j} f S_{j} g. 
\end{equation}

%%%%%%%%%%%%%%%%%%%%%
\subsection{Compactness}
%%%%%%%%%%%%%%%%%%%%
Since we look for weak solutions, we use compactness arguments when we pass to the limit in weak formulations. 

\begin{lemma}\cite{Temam} \label{lemma:2.1}
Let $X_0,X,X_1$ be reflexive Banach spaces such that
\[
X_0\subset \subset X\subset X_1,
\] 
where $X_0$ is compactly embedded in $X$. Let $T>0$ be a finite number and let $\alpha_0$ and $\alpha_1$ be two finite numbers such that $\alpha_i>1$. Then, $\displaystyle Y=\left\{u\in L^{\alpha_0}\left(0,T; X_0\right),\  \partial_t u\in L^{\alpha_1}\left(0,T; X_1\right)\right\}$ is compactly embedded in $L^{\alpha_0}\left(0,T; X\right)$.
\end{lemma}

\begin{lemma}[\cite{constantin2010global}]\label{IVlemapaco}
Consider a sequence $(\theta^{\epsilon}) \in C([0,T]\times B_{R}(0))$ that is uniformly bounded in $L^\infty([0,T],W^{1,\infty}(B_{R}(0)))$. Assume further that the weak derivative $\frac{d\theta^{\epsilon}}{dt}$ is in $L^\infty([0,T],L^\infty(B_{R}(0)))$ (not necessarily uniform) and is uniformly bounded in $L^\infty([0,T],W^{-2,\infty}_*(B_{R}(0)))$. Finally suppose that $\theta^{\epsilon}_x\in C([0,T]\times B_{R}(0))$. Then there exists a subsequence of $(\theta^{\epsilon})$ that converges strongly in $L^\infty([0,T]\times B_{R}(0))$.
\end{lemma}

%%%%%%%%%%%%%%%%%%%%%%%%
\section{Proof of Theorem \ref{GW 1}}
%%%%%%%%%%%%%%%%%%%%%%%%

%%%%%%%%%%%%%%%%%
\subsection{A priori estimates}
%%%%%%%%%%%%%%%%
We first obtain a priori bounds of the equation
\eqn\label{Hilbert}
\theta_t+\left(\mathcal{H}\theta\right)\theta_x +\delta\theta \Lambda\theta+\Lambda^{\gamma}\theta=0, 
\een
We note that by the minimum principle applied to (\ref{Hilbert}), we have $\theta(t,x)\ge 0$ for all $t\ge 0$.

To obtain $H^{\frac{1}{2}}$ bound of $\theta$, we begin with the $L^{2}$ bound. We multiply (\ref{Hilbert}) by $\theta$ and integrate over $\mathbb{R}$. Then, 
\[
\frac{1}{2}\frac{d}{dt}\left\| \theta\right\|^{2}_{L^{2}}+\left\|\Lambda^{\frac{\gamma}{2}}\theta\right\|^{2}_{L^{2}} =-\int \left[\left(\mathcal{H}\theta\right)\theta_{x} \theta\right] dx -\delta\int \left[\theta^{2} \Lambda \theta\right]dx =\left(\frac{1}{2}-\delta\right)\int \left[\theta^{2} \Lambda \theta\right]dx.
\]
Since $\theta\ge 0$, we have  
\[
\int \left[\theta^{2} \Lambda \theta\right]dx=\int  \int \frac{\left(\theta(x)-\theta(y)\right)^{2}}{|x-y|^{2}} \cdot \frac{\theta(x)+\theta(y)}{2} dxdy\ge 0
\]
and thus  
\eqn \label{L2 bound of theta}
\left\|\theta(t) \right\|^{2}_{L^{2}}+2\int^{t}_{0}\left\|\Lambda^{\frac{\gamma}{2}}\theta(s)\right\|^{2}_{L^{2}}ds\leq \left\| \theta_{0} \right\|^{2}_{L^{2}}.
\een

We next estimate $\theta$ in $\dot{H}^{\frac{1}{2}}$. We multiply (\ref{Hilbert}) by $\Lambda \theta$ and integrate over $\mathbb{R}$:
\[
\frac{1}{2}\frac{d}{dt}\left\| \Lambda^{\frac{1}{2}}\theta \right\|^{2}_{L^{2}}+\left\|\Lambda^{\frac{1+\gamma}{2}}\theta\right\|^{2}_{L^{2}}=-\int \left[\left(\mathcal{H}\theta\right)\theta_{x}\Lambda \theta\right]dx -\delta\int \left[\theta \left(\Lambda \theta\right)^{2}\right]dx.
\]
By (\ref{Hilbert identity}), we have
\begin{equation*} 
 \begin{split}
 -\int \left[\left(\mathcal{H}\theta\right)\theta_{x}\Lambda \theta\right]dx = \int \left[\theta \mathcal{H} \left(\theta_{x}\left(\mathcal{H}\theta_{x}\right)\right)\right]dx =\frac{1}{2}\int \left[\theta\left(\left(\Lambda \theta\right)^{2}-\left(\theta_{x}\right)^{2}\right)\right]dx,
  \end{split}
\end{equation*}
and hence 
\begin{equation*}
\begin{split}
&\frac{1}{2}\frac{d}{dt}\left\| \Lambda^{\frac{1}{2}}\theta \right\|^{2}_{L^{2}}+\left\|\Lambda^{\frac{1+\gamma}{2}}\theta\right\|^{2}_{L^{2}}=\left(\frac{1}{2}-\delta\right)\int \left[\theta\left(\Lambda \theta \right)^{2}\right]dx- \frac{1}{2}\int \left[\theta \left(\theta_{x}\right)^{2}\right]dx\leq 0,
\end{split}
\end{equation*}
where we use the sign conditions $\theta\ge 0$ and $\delta\ge \frac{1}{2}$. This leads to the inequality   
\eqn \label{H1/2 bound of theta}
\left\| \Lambda^{\frac{1}{2}}\theta(t) \right\|^{2}_{L^{2}}+2\int^{t}_{0}\left\|\Lambda^{\frac{1+\gamma}{2}}\theta(s)\right\|^{2}_{L^{2}}ds\leq \left\| \Lambda^{\frac{1}{2}}\theta_{0} \right\|^{2}_{L^{2}}.
\een

By (\ref{L2 bound of theta}) and (\ref{H1/2 bound of theta}), we obtain that 
\eqn \label{H1/2 bound of theta final}
\left\| \theta(t) \right\|^{2}_{H^{\frac{1}{2}}}+2\int^{t}_{0}\left\|\Lambda^{\frac{\gamma}{2}}\theta(s)\right\|^{2}_{H^{\frac{1}{2}}}ds\leq \left\|\theta_{0} \right\|^{2}_{H^{\frac{1}{2}}}.
\een

We finally estimate $\theta$ in $L^{1}$. Since $\theta\ge 0$,
\[
\frac{d}{dt} \|\theta\|_{L^{1}}=\frac{d}{dt}\int \theta dx=(1-\delta)\int \theta \Lambda\theta dx \leq C \left\|\Lambda^{\frac{1}{2}}\theta\right\|^{2}_{L^{2}}
\]
and thus we conclude that 
\eqn \label{L1 bound of theta}
\|\theta(t)\|_{L^{1}}\leq \|\theta_{0}\|_{L^{1}}+ C \int^{t}_{0}\left\| \theta(s) \right\|^{2}_{\dot{H}^{\frac{1}{2}}}ds \leq \|\theta_{0}\|_{L^{1}}+ Ct\left\|\theta_{0} \right\|^{2}_{H^{\frac{1}{2}}}.
\een

%%%%%%%%%%%%%%%%%%%%%%%%%%
\subsection{Approximation and passing to limit}
%%%%%%%%%%%%%%%%%%%%%%%%%%
We first regularize initial data as $\theta^{\epsilon}_{0}=\rho_{\epsilon}\ast \theta_{0}$ where $\rho_{\epsilon}$ is a standard mollifier. We then regularize the equation by putting the Laplacian with the coefficient $\epsilon$:
\eqn \label{regularized equation}
\theta^{\epsilon}_t+\left(\mathcal{H}\theta^{\epsilon}\right)\theta^{\epsilon}_x + \delta\theta^{\epsilon} \Lambda\theta^{\epsilon}+\Lambda^{\gamma}\theta^{\epsilon}=\epsilon \theta^{\epsilon}_{xx}.
\een
For the proof of the existence of a global-in-time smooth solution, see \cite{Lazar} (Section 6). Moreover, $(\theta^{\epsilon})$ satisfies that
\[
\left\| \theta^{\epsilon}(t) \right\|_{L^{1}}+ \left\| \theta^{\epsilon}(t) \right\|^{2}_{H^{\frac{1}{2}}}+2\int^{t}_{0}\left\|\Lambda^{\frac{\gamma}{2}}\theta^{\epsilon}(s)\right\|^{2}_{H^{\frac{1}{2}}}ds +\epsilon\left\|\nabla \theta^{\epsilon}\right\|^{2}_{H^{\frac{1}{2}}} \leq \left\|\theta_{0}\right\|_{L^{1}}+ C(1+t)\left\|\theta_{0} \right\|^{2}_{H^{\frac{1}{2}}}.
\]
Therefore, $(\theta_{\epsilon})$ is bounded in $\mathcal{A}_{T}$ uniformly in $\epsilon>0$. From this, we have uniform bounds  
\[
\mathcal{H}\theta^{\epsilon}\in L^{4}\left(0,T; L^{4}\right), \quad \theta^{\epsilon}\in L^{2}\left(0,T; L^{4}\right)
\]
and hence
\[
\left(\left(\mathcal{H}\theta^{\epsilon}\right)\theta^{\epsilon}\right)_x\in L^{\frac{4}{3}}\left(0,T; H^{-1}\right).
\]
Moreover, 
\[
\Lambda^{\gamma}\theta^{\epsilon} +\epsilon \theta^{\epsilon}_{xx}\in L^{2}\left(0,T; H^{-1}\right).
\]
To estimate $\theta^{\epsilon} \Lambda\theta^{\epsilon}$, we use the duality argument. For any $\chi \in L^{2}\left(0,T; H^{2}\right)$,
\begin{equation*}
\begin{split}
\left|\langle \theta^{\epsilon} \Lambda \theta^{\epsilon}, \chi \rangle\right| &\leq \int \left|\widehat{\theta^{\epsilon}\Lambda\theta^{\epsilon}}(\xi) \widehat{\chi}(\xi)\right|d\xi \leq \int \int \left|\widehat{\theta^{\epsilon}}(\xi-\eta)\right| |\eta| \left|\widehat{\theta^{\epsilon}}(\eta)\right| \left|\widehat{\chi}(\xi)\right|d\eta d\xi\\
& \leq \int \int \left|\widehat{\theta^{\epsilon}}(\xi-\eta)\right| \left(|\eta|^{\frac{1}{2}}\left(|\xi-\eta|^{\frac{1}{2}}+|\xi|^{\frac{1}{2}}\right)\right) \left|\widehat{\theta^{\epsilon}}(\eta)\right| \left|\widehat{\chi}(\xi)\right|d\eta d\xi\\
& =\int \int |\xi-\eta|^{\frac{1}{2}}\left|\widehat{\theta^{\epsilon}}(\xi-\eta)\right| |\eta|^{\frac{1}{2}} \left|\widehat{\theta^{\epsilon}}(\eta)\right| \left|\widehat{\chi}(\xi)\right|d\eta d\xi \\
&+ \int \int \left|\widehat{\theta^{\epsilon}}(\xi-\eta)\right| |\eta|^{\frac{1}{2}} \left|\widehat{\theta^{\epsilon}}(\eta)\right| |\xi|^{\frac{1}{2}}\left|\widehat{\chi}(\xi)\right|d\eta d\xi\\
&\leq \left\|\Lambda^{\frac{1}{2}}\theta^{\epsilon}\right\|^{2}_{L^{2}}\left\|\widehat{\chi}\right\|_{L^{1}}+ \left\|\Lambda^{\frac{1}{2}}\theta^{\epsilon}\right\|_{L^{2}}\left\|\widehat{\theta}\right\|_{L^{1}} \left\|\Lambda^{\frac{1}{2}}\chi\right\|_{L^{2}}  \leq C\left\|(1+\Lambda)^{\frac{1}{2}+\frac{\gamma}{4}}\theta^{\epsilon}\right\|^{2}_{L^{2}}\|\chi\|_{H^{2}}.
\end{split}
\end{equation*}
Since $(1+\Lambda)^{\frac{1}{2}+\frac{\gamma}{4}}\theta^{\epsilon} \in L^{4}\left(0,T; L^{2}\right)$ uniformly in $\epsilon>0$, we have 
\[
\int \left|\langle \theta^{\epsilon} \Lambda \theta^{\epsilon}, \chi \rangle\right|dt\leq C \left\|(1+\Lambda)^{\frac{1}{2}+\frac{\gamma}{4}}\theta^{\epsilon}\right\|^{2}_{L^{4}_{T}L^{2}} \|\chi\|_{L^{2}_{T}H^{2}}.
\]
This implies that $\theta^{\epsilon}\Lambda \theta^{\epsilon} \in L^{2}\left(0,T; H^{-2}\right)$. So, we conclude that  from the equation of $\theta^{\epsilon}_{t}$
\[
\theta^{\epsilon}_{t}\in L^{\frac{4}{3}}\left(0,T; H^{-2}\right).
\]

We now extract a subsequence of $\left(\theta^{\epsilon}\right)$, using the same index $\epsilon$ for simplicity, and a function $\theta \in \mathcal{A}_T$ such that 
\begin{equation}\label{convergence of approximate solutions}
\begin{split}
& \theta^{\epsilon} \stackrel{\star}{\rightharpoonup} \theta \quad \text{in} \quad L^{\infty}\left(0,T; L^{p}\cap H^{\frac{1}{2}}\right) \quad \text{for all $p\in (1,\infty)$},\\
& \theta^{\epsilon} \rightharpoonup \theta \quad \text{in} \quad L^{2}\left(0,T; H^{\frac{\gamma+1}{2}}\right),\\
& \theta^{\epsilon} \rightarrow  \theta \quad \text{in} \quad L^{2}\left(0,T; H^{\frac{1}{2}}\right),\\
& \theta^{\epsilon}\rightarrow \theta \quad \text{in} \quad L^{2}\left(0,T; L^{p}_{\text{loc}}\right) \quad \text{for all $p\in (1,\infty)$}
\end{split}
\end{equation}
where we use Lemma \ref{lemma:2.1} for the strong convergence.

\textcolor{black}{
We now multiply (\ref{regularized equation}) by a test function $\psi\in \mathcal{C}^{\infty}_{c}\left([0,T)\times\mathbb{R}\right)$ and integrate over $\mathbb{R}$. Then,
\[
\int^{T}_{0}\int \left[\theta^{\epsilon} \psi_{t} + \left(\mathcal{H}\theta^{\epsilon}\right)\theta^{\epsilon}\psi_{x}+(1-\delta)\Lambda\theta^{\epsilon}\theta^{\epsilon} \psi -\theta^{\epsilon}\Lambda^{\textcolor{red}{\gamma}}\psi+ \epsilon \theta^{\epsilon}\psi_{xx}\right]  dxdt =\int\theta^{\epsilon}_{0}(x)\psi(0,x)dx\,,
\]
which can be rewritten as
\begin{equation} \label{weak form with commutator}
\begin{split}
&\int^{T}_{0}\int \Big[\theta^{\epsilon} \psi_{t} +  \underbrace{\left(\mathcal{H}\theta^{\epsilon}\right)\theta^{\epsilon}\psi_{x}}_{\text{I}} -\theta^{\epsilon} \Lambda^{\textcolor{red}{\gamma}}\psi+ \epsilon \theta^{\epsilon}\psi_{xx}\Big]  dxdt - \int \theta^{\epsilon}_{0}(x)\psi(0,x)dx\\
& = -(1-\delta)\int^{T}_{0}\int \underbrace{\Lambda^{\frac{1}{2}}\theta^{\epsilon} \left[\Lambda^{\frac{1}{2}},\psi \right]\theta^{\epsilon}}_{\text{II}}  dxdt -(1-\delta) \int^{T}_{0}\int \underbrace{\left|\Lambda^{\frac{1}{2}}\theta^{\epsilon}\right|^{2}\psi}_{\text{III}}  dxdt.
\end{split}
\end{equation}
}
By Lemma \ref{lemma:2.1} with
\[
X_0=L^2\left(0,T; H^{\frac{1}{2}}\right),\quad X=L^2\left(0,T; L^2_{\text{loc}}\right),\quad X_1=L^2\left(0,T; H^{-2}\right),
\]
we can pass to the limit to $\text{I}$. Moreover, since
\[
\left[\Lambda^{\frac{1}{2}},\psi \right]\theta^{\epsilon} \rightarrow \left[\Lambda^{\frac{1}{2}},\psi \right]\theta
\]
strongly in $L^{2}\left(0,T; L^{6}\right)$ by Lemma \ref{commutator lemma} and $\Lambda^{\frac{1}{2}}\theta^{\epsilon} $ converges weakly in $L^{2}\left(0,T; L^{2}\right)$ by (\ref{convergence of approximate solutions}), we can pass to the limit to $\text{II}$. Lastly, Lemma \ref{lemma:2.1} with
\[
X_0=L^2\left(0,T; H^{\frac{1+\gamma}{2}}\right),\quad X=L^2\left(0,T; H^{\frac{1}{2}}_{\text{loc}}\right),\quad X_1=L^2\left(0,T; H^{-2}\right),
\]
allows to pass to the limit  to $\text{III}$. Combining all the limits together, we obtain that 
\begin{equation} \label{weak form with commutator 2}
\begin{split}
\int^{T}_{0}\int \left[\theta \psi_{t} +\left(\mathcal{H}\theta\right)\theta\psi_{x}+(1-\delta)\Lambda\theta\theta^{\epsilon} \psi\right]  dxdt =\int \theta_{0}(x)\psi(0,x)dx.
\end{split}
\end{equation}

%%%%%%%%%%%%%%%%%%%%%%%%%
\subsection{Uniqueness when $\gamma=1$}
%%%%%%%%%%%%%%%%%%%%%%%%%
To show the uniqueness of a weak solution, let $\theta=\theta_{1}-\theta_{2}$. Then, $\theta$ satisfies the following equation:
\eqn \label{eq:2.7}
\theta_{t}+\Lambda \theta=-\left(\mathcal{H}\theta\right) \theta_{1x}- \left(\mathcal{H}\theta_{2}\right) \theta_{x}-\delta \theta \Lambda \theta_{1}-\delta \theta_{2}\Lambda \theta,  \quad \theta(0,x)=0.  
\een
We multiply $\theta$ to (\ref{eq:2.7}) and integrate over $\mathbb{R}$. Then, 
\begin{equation*} 
 \begin{split}
 \frac{1}{2}\frac{d}{dt}\left\|\theta\right\|^{2}_{L^{2}}+\left\|\Lambda^{\frac{1}{2}}\theta \right\|^{2}_{L^{2}}  = \int \left[-\left(\mathcal{H}\theta\right) \theta_{1x}- \left(\mathcal{H}\theta_{2}\right) \theta_{x}-\delta \theta \Lambda \theta_{1}-\delta \theta_{2}\Lambda \theta\right]\theta dx.
 \end{split}
\end{equation*}
The first three terms in the right-hand side are easily bounded by
\[
C\left\|\theta_{1x}\right\|_{L^{2}} \|\theta\|^{2}_{L^{4}}+ C\left\|\theta_{2x}\right\|_{L^{2}} \|\theta\|^{2}_{L^{4}}.
\]
Moreover, the last term is bounded by using Lemma \ref{pointwise lemma}
\[
-\delta\int  \theta_{2}\theta\Lambda \theta dx\leq -\frac{\delta}{2} \int \theta_{2}\Lambda \theta^{2} dx= -\frac{\delta}{2} \int \theta^{2}\Lambda \theta_{2} dx\leq C \left\|\theta_{2x}\right\|_{L^{2}} \|\theta\|^{2}_{L^{4}}.
\]
Hence we derive that 
\begin{equation*} 
 \begin{split}
\frac{d}{dt}\left\|\theta\right\|^{2}_{L^{2}}+\left\|\Lambda^{\frac{1}{2}}\theta \right\|^{2}_{L^{2}}  &\leq C\left\|\theta_{1x}\right\|_{L^{2}} \|\theta\|^{2}_{L^{4}}+ C\left\|\theta_{2x}\right\|_{L^{2}} \|\theta\|^{2}_{L^{4}}\\
 & \leq C \left(\left\|\theta_{1x}\right\|^{2}_{L^{2}}+ \left\|\theta_{2x}\right\|^{2}_{L^{2}} \right) \|\theta\|^{2}_{L^{2}} +\frac{1}{2}\left\|\Lambda^{\frac{1}{2}}\theta \right\|^{2}_{L^{2}}.
 \end{split}
\end{equation*}
Since 
\[
\theta_{1x} \in L^{2}\left(0,T: L^{2}\right), \quad \theta_{2x} \in L^{2}\left(0,T: L^{2}\right)
\]
when $\gamma=1$, we conclude that $\theta=0$ in $L^{2}$ and thus a weak solution is unique. This completes the proof of Theorem \ref{GW 1}.

%%%%%%%%%%%%%%%%%%%%%%%%%%%%%
\section{Proof of Theorem \ref{GW infinite energy}}
%%%%%%%%%%%%%%%%%%%%%%%%%%%%%

%%%%%%%%%%%%%%%%
\subsection{A priori estimate}
%%%%%%%%%%%%%%%%
We consider the equation
\eqn \label{equation infinite 1}
\theta_{t}+\left(\mathcal{H}\theta\right)\theta_{x}+\delta \theta \Lambda \theta+\Lambda \theta=0.
\een
Since (\ref{equation infinite 1}) satisfies the minimum and maximum principles, we have 
\[
\theta(t,x)\ge 0, \quad \|\theta(t)\|_{L^\infty}\leq \|\theta_0\|_{L^\infty}.
\]
We begin with the $L^{2}(w_{\beta}dx)$ bound. For notational simplicity, we suppress the dependence on $\beta$. We multiply (\ref{equation infinite 1}) by $\theta w$ and integrate in $x$. Then, 
\begin{equation*}
\begin{split}
\frac{1}{2}\frac{d}{dt} \left\|\theta\right\|^{2}_{L^{2}(wdx)} &+ \left\|\Lambda^{1/2}\theta\right\|^{2}_{L^{2}(wdx)}=-\int \left(\mathcal{H}\theta\right)\theta_{x} \theta wdx -\delta\int \theta \left(\Lambda \theta\right) \theta w dx-\int \Lambda^{1/2}\theta \left[\Lambda^{1/2},w\right]\theta dx\\
& =-\frac{1}{2} \left(\mathcal{H}\theta\right)\left(\theta^{2}\right)_{x} wdx -\delta\int \theta^{2} \left(\Lambda \theta\right) w dx -\int \Lambda^{1/2}\theta \left[\Lambda^{1/2},w\right]\theta dx\\
&=\left(\frac{1}{2}-\delta\right)\int \theta^{2} \left(\Lambda \theta\right) w dx+\frac{1}{2}\int \left(\mathcal{H}\theta\right) \theta^{2}w_{x}dx-\int \Lambda^{1/2}\theta \left[\Lambda^{1/2},w\right]\theta dx \\
& \leq C( \Vert \theta_{0} \Vert_{L^{\infty}}, \delta) \left(\Vert \sqrt{w}\theta \Vert_{L^{2}} \Vert \sqrt{w}\Lambda\theta \Vert_{L^{2}} + \Vert \sqrt{w} \theta \Vert^{2}_{L^{2}} \right) + \int \sqrt{w} \vert \theta \vert \ \frac{1}{\sqrt{w}}\left\vert \left[\Lambda^{1/2},w\right]\theta \right \vert \ dx 
\end{split}
\end{equation*}
Then, using lemma \ref{lem24} one obtains, for $C( \Vert \theta_{0} \Vert_{L^{\infty}}, \delta)=C_1=(1+\delta)\Vert \theta_{0} \Vert_{L^{\infty}}$ for all $\eta_1>0$ and 
\begin{eqnarray*}
\frac{1}{2}\frac{d}{dt} \left\|\theta\right\|^{2}_{L^{2}(wdx)} + \left\|\Lambda^{\frac{1}{2}}\theta\right\|^{2}_{L^{2}(wdx)} 
&\leq&  \left(\Vert\theta \Vert_{L^{2}(\w)} \Vert \Lambda\theta \Vert_{L^{2}_{w}} + \Vert \theta \Vert^{2}_{L^{2}_{w}} \right)C_1  + \Vert \theta \Vert^{2}_{L^{2}(\w)}  \\
&\leq&  \left(\frac{C_1}{2\eta_1}+C_1+1\right)  \Vert\theta \Vert^{2}_{L^{2}(\w)} + \frac{C_1\eta_1}{2}\Vert \Lambda\theta \Vert^{2}_{L^{2}_{w}} 
\end{eqnarray*}
In particular,
\begin{eqnarray} \label{L2b}
\frac{1}{2}\frac{d}{dt} \left\|\theta\right\|^{2}_{L^{2}(wdx)}&\leq& \left(\frac{C_1}{2\eta_1}+C_1+1\right)  \Vert\theta \Vert^{2}_{L^{2}(\w)} + \frac{C_1\eta_1}{2}\Vert \Lambda\theta \Vert^{2}_{L^{2}_{w}} 
\end{eqnarray}

We next multiply (\ref{equation infinite 1}) by $\Lambda^{\frac{1}{2}}\left(w\Lambda^{\frac{1}{2}}\theta\right)$ and integrate in $x$.  We first focus one the term $\delta \theta \Lambda \theta$, the resulting term in the computation of the evolution of the $\dot H^{1/2}(\w)$ norm  is
\begin{eqnarray*}
T_{\delta} \theta&=&\delta\int \Lambda^{1/2}(w \Lambda^{1/2} \theta) \theta \Lambda \theta \ dx \\
&=& \delta \int \theta \sqrt{w}\Lambda \theta \ \frac{1}{\sqrt{w}}\left[\Lambda^{1/2}, w \right] \Lambda^{1/2}\theta \ dx + \int  \theta \ w\vert \Lambda \theta \vert^{2} \ dx \\
&\leq&  \delta \Vert \theta_{0} \Vert_{L^{\infty}} \left(\Vert \theta \Vert_{\dot H^{1}(\w)} \Vert \theta \Vert_{\dot H^{1/2}(\w)} + \Vert \theta \Vert^{2}_{\dot H^{1}(\w)} \right) \\
\end{eqnarray*}
hence 
\begin{eqnarray}\label{h1}
\vert T_{\delta} \theta \vert \leq \label{h1} \left(\frac{\eta_2}{2}+1\right) C_2\Vert \theta \Vert^{2}_{\dot H^{1}(\w)} +\frac{2C_2}{\eta_2} \Vert \theta \Vert^{2}_{\dot H^{1/2}(\w)},
\end{eqnarray}
where $C_2=\delta\Vert \theta_{0} \Vert_{L^{\infty}}$. \\

\noindent For the transport part (which is the case $\delta=0$), we only need to use inequality 4.4 of \cite{Lazar}, indeed, it is shown that, for some constants $C_3>0,C_4>0$ one has
\begin{eqnarray}\label{eqsob2}
 \frac{1}{2}\frac{d}{dt}  \Vert \theta \Vert^{2}_{\dot H^{1/2}(\w)} \ dx  &\leq& \nonumber (C_3 \Vert \theta_{0} \Vert_{\infty}-1)\Vert \theta \Vert^{2}_{\dot H^{1}(\w)} \\
 &+&  C_4 \left(\int \theta^2 \ w\ dx + \int \vert\Lambda^{1/2}\theta\vert^2 \ \ w \ dx\right).
\end{eqnarray}
Hence, we get the following control for the full equation \ref{equation infinite 1}  by collecting the estimates \ref{L2b}, \ref{h1} and \ref{eqsob2}
\begin{equation}\label{eqsob2}\begin{split}
\frac{1}{2}\frac{d}{dt} \Vert \theta \Vert^{2}_{H^{1/2}(\w)} \leq& \left(C_3 \Vert \theta_{0} \Vert_{\infty}+ \left(\frac{\eta_2}{2}+1\right) C_1 + \frac{C_1\eta_1}{2}-1\right)\Vert \theta \Vert^{2}_{\dot H^{1}(\w)} \\
 &+ (C_4 +\frac{2C_2}{\eta_2}) \Vert \theta \Vert^{2}_{\dot H^{1/2}(\w)} 
 \end{split}
\end{equation}
 Hence, choosing $\eta_1, \eta_2$ and $\Vert \theta_{0} \Vert_{\infty}$ small enough (for instance less than $\frac{1}{1+\delta} \frac{1}{100}$ so that $C_1$ would be also small), one gets
 \begin{equation}
\begin{split}
\frac{1}{2}\frac{d}{dt}  \Vert \theta \Vert^2_{H^{1/2}(\w)} \ dx + C_6 \Vert \theta \Vert^{2}_{\dot H^{1}(w)}  \leq&  C_7\Vert \theta \Vert^{2}_{\dot H^{1/2}(\w)}.
 \end{split}
 \end{equation}
Hence, integrating in time $s\in[0,T]$ and using  Gronwall's inequality  one obtains that for all finite $T>0$ 
\begin{equation}\label{L2 and H1/2 bound of theta w}
\begin{split}
\left\|\theta(T)\right\|^{2}_{H^{1/2}(wdx)}+ \Vert \theta \Vert^{2}_{L^{2}([0,T],\dot H^{1}(\w))}\leq \left\|\theta_{0}\right\|^{2}_{H^{1/2}(wdx)}\exp\left(C_{2} \Vert\theta_{0}\Vert_{L^{\infty}}T\right).
\end{split}
\end{equation}

%%%%%%%%%%%%%%%%%%%%%%%%%%
\subsection{Approximation and passing to limit}
%%%%%%%%%%%%%%%%%%%%%%%%%%
To show the existence of a weak solution in $\mathcal{D}_{T}$, we first approximate the initial data $\theta_{0}$. Let $\chi$ be a smooth positive function such that $\chi(x)=1$ for $|x|\leq 1$ and $\chi(x)=0$ for $|x|\ge 2$. Let $\chi_{R}(x)=\chi(x/N)$, $N\in \mathbb{N}$, and consider truncated initial data $\theta^{N}_{0}(x)=\theta_{0} (x)\chi_{N}(x)$. Then, a direct computation shows that
\[
\lim_{N\rightarrow \infty}\left\|\theta^{N}_{0}-\theta_{0}\right\|_{H^{\frac{1}{2}}(wdx)}=0.
\]
Moreover, this truncation does not alter the non-negativity and does not increase the $L^{\infty}$ norm. So, if $\|\theta_{0}\|_{L^{\infty}}$ is sufficiently small, there is a global-in-time solution of
\begin{equation}
\begin{split}
& \partial_{t}\theta^{N}+\mathcal{H}\theta^{N} \partial_{x}\theta^{N}+\delta \theta^{N}\Lambda \theta^{N}+ \Lambda \theta^{N}=0,\quad \theta^{N}(0,x)=\theta^{N}_{0}(x).
\end{split}
\end{equation}
From the a priori estimates, the sequence $(\theta^{N})$ is bounded in  
\[
L^{\infty}([0,T], H^{\frac{1}{2}}(wdx)) \cap L^{2}([0,T], H^{1}(wdx))
\]
uniformly with respect to $N$. We now take a test function $\psi\in \mathcal{C}^{\infty}_{c}\left([0,T)\times\mathbb{R}\right)$. Then, $\psi\theta^{N}$ is bounded in $L^{2}([0,T], H^{1})$. Moreover, since $\theta^{N}\in L^{\infty}([0,T]\times \mathbb{R})$
\begin{equation*}
\begin{split}
&\psi \left(\mathcal{H}\theta^{N} \partial_{x}\theta^{N} +\delta \theta^{N}\Lambda \theta^{N}+ \Lambda \theta^{N}\right)\\
&=(\psi \mathcal{H}\theta^{N} \theta^{N})_{x}-\psi_{x}\mathcal{H}\theta^{N} \theta^{N} +(\delta-1)\psi \theta^{N}\Lambda\theta^{N} +\psi \Lambda \theta^{N} \in L^{2}([0,T], H^{-1}).
\end{split}
\end{equation*}
By Lemma \ref{lemma:2.1}, we can pass to the limit to the weak formulation, 
\[
\int^{T}_{0}\int \left[\theta^{N} \psi_{t} + \left(\mathcal{H}\theta^{N}\right)\theta^{N}\psi_{x}+(1-\delta)\Lambda\theta^{N}\theta^{N} \psi -\theta^{N}\Lambda\psi\right]  dxdt =\int\theta^{N}_{0}(x)\psi(0,x)dx,
\]
to obtain a weak solution $\theta$ which is also in 
\[
L^{\infty}([0,T], H^{\frac{1}{2}}(wdx)) \cap L^{2}([0,T], H^{1}(wdx)).
\]

%%%%%%%%%%%%%%
\subsection{Uniqueness}
%%%%%%%%%%%%%%
To show the uniqueness of a weak solution, we consider the equation of $\theta=\theta_{1}-\theta_{2}$ given by 
\eqn \label{eq:2.77}
\theta_{t}+\Lambda \theta=-\left(\mathcal{H}\theta\right) \theta_{1x}- \left(\mathcal{H}\theta_{2}\right) \theta_{x}-\delta \theta \Lambda \theta_{1}-\delta \theta_{2}\Lambda \theta,  \quad \theta(0,x)=0.  
\een
We multiply $w\theta$ to (\ref{eq:2.77}) and integrate over $\mathbb{R}$. Then, 
\begin{equation*} 
 \begin{split}
 \frac{1}{2}\frac{d}{dt}\left\|\theta\right\|^{2}_{L^{2}(wdx)}+\left\|\Lambda^{\frac{1}{2}}\theta \right\|^{2}_{L^{2}(wdx)}  &= \int \left[-\left(\mathcal{H}\theta\right) \theta_{1x}- \left(\mathcal{H}\theta_{2}\right) \theta_{x}-\delta \theta \Lambda \theta_{1}-\delta \theta_{2}\Lambda \theta\right]\theta wdx\\
 & -\int \Lambda^{\frac{1}{2}}\theta \left[\Lambda^{\frac{1}{2}},w\right]\theta dx.
 \end{split}
\end{equation*}
As before, the last term is bounded by
\[
\int \Lambda^{\frac{1}{2}}\theta \left[\Lambda^{\frac{1}{2}},w\right]\theta dx \leq \frac{1}{2} \left\|\Lambda^{\frac{1}{2}}\theta\right\|^{2}_{L^{2}(wdx)} +  C \left\|\theta\right\|^{2}_{L^{2}(wdx)}.
\]
The first three terms in the right-hand side are easily bounded by
\[
C\left(\left\|\theta_{1x}\right\|_{L^{2}(wdx)}+\left\|\theta_{2x}\right\|_{L^{2}(wdx)}+\left\|\theta_{2}\right\|_{L^{2}(wdx)} \right) \|\theta\|^{2}_{L^{4}(wdx)}.
\]
Moreover, since $\delta>0$, $\theta_{2}\ge 0$ and $w\ge 0$, the fourth term is bounded by using Lemma \ref{pointwise lemma}
\begin{equation*}
\begin{split}
-\delta\int  \theta_{2}\theta\Lambda \theta wdx &\leq -\frac{\delta}{2} \int \theta_{2} w\Lambda \theta^{2} dx= -\frac{\delta}{2} \int \theta^{2}\Lambda (\theta_{2}w) dx = \frac{\delta}{2} \int \mathcal{H}(\theta^{2})(\theta_{2}w)_{x} dx \\
& \leq C \left(\left\|\theta_{2x}\right\|_{L^{2}(wdx)}+ \left\|\theta_{2}\right\|_{L^{2}(wdx)}\right) \|\theta\|^{2}_{L^{4}(wdx)}.
\end{split}
\end{equation*}
Hence we obtain that 
\begin{equation*} 
 \begin{split}
& \frac{d}{dt}\left\|\theta\right\|^{2}_{L^{2}(wdx)}+\left\|\Lambda^{\frac{1}{2}}\theta \right\|^{2}_{L^{2}(wdx)} \\
&\leq C\left(\left\|\theta_{1x}\right\|_{L^{2}(wdx)}+\left\|\theta_{2x}\right\|_{L^{2}(wdx)}+\left\|\theta_{2}\right\|_{L^{2}(wdx)} \right) \|\theta\|^{2}_{L^{4}(wdx)} +  C \left\|\theta\right\|^{2}_{L^{2}(wdx)}\\
 & \leq C \left(1+ \left\|\theta_{1x}\right\|^{2}_{L^{2}(wdx)}+ \left\|\theta_{2x}\right\|^{2}_{L^{2}(wdx)} + \left\|\theta_{2}\right\|^{2}_{L^{2}(wdx)}\right) \|\theta\|^{2}_{L^{2}(wdx)} +\frac{1}{2}\left\|\Lambda^{\frac{1}{2}}\theta \right\|^{2}_{L^{2}(wdx)},
 \end{split}
\end{equation*}
where we use (\ref{GN weight}) to obtain the last inequality. Since 
\[
\theta_{1x} \in L^{2}\left(0,T: L^{2}(wdx)\right), \quad \theta_{2x} \in L^{2}\left(0,T: H^{1}(wdx)\right),
\]
we conclude that $\theta=0$ in $L^{2}(wdx)$ and thus a weak solution is unique. This completes the proof of Theorem \ref{GW infinite energy}. \\

In the two next subsections we shall prove a global existence theorem and a local existence theorem, for respectively the subcritical case and the supercritical case), we shall just focus on the {\it{a priori}} estimates since the construction by compactness is classical (see \cite{Lazar2}).

%%%%%%%%%%%%%%%%%%%%%%%%%%%%%
\section{Proof of Theorem \ref{sub}}
%%%%%%%%%%%%%%%%%%%%%%%%%%%%%

%%%%%%%%%%%%%%%%
\subsection{A priori estimate}
%%%%%%%%%%%%%%%%

Since the case $\delta=0$ has been already done in \cite{Lazar2}, one observes that it suffices to estimate the $L^{2}(\w)$ part 
$$
T_1=\delta \int w_{\beta,k}   \theta^{2} \Lambda \theta \ dx,
$$
and the weighted homogeneous Sobolev part, that is
$$
T_2=\delta \int w_{\beta,k}  \Lambda^{s} \theta \ \Lambda^{s}(\theta \Lambda \theta) \ dx.
$$
We shall control $T_1$ and $T_2$, for the sake of readibility, we shall just write $w$ in stead of $w_{\beta,k}$. We first control the $L^{2}$ part
\begin{equation*}
\begin{split}
\frac{1}{2}\frac{d}{dt} \left\|\theta\right\|^{2}_{L^{2}(wdx)} &+ \left\|\Lambda^{\alpha/2}\theta\right\|^{2}_{L^{2}(wdx)}=-\int \left(\mathcal{H}\theta\right)\theta_{x} \theta wdx -\delta\int \theta \left(\Lambda \theta\right) \theta w dx-\int \Lambda^{\alpha/2}\theta \left[\Lambda^{\alpha/2},w\right]\theta dx\\
& =-\frac{1}{2}\int \left(\mathcal{H}\theta\right)\left(\theta^{2}\right)_{x} wdx -\delta\int \theta^{2} \left(\Lambda \theta\right) w dx -\int \Lambda^{\alpha/2}\theta \left[\Lambda^{\alpha/2},w\right]\theta dx\\
&=\left(\frac{1}{2}-\delta\right)\int \theta^{2} \left(\Lambda \theta\right) w dx+\frac{1}{2}\int \left(\mathcal{H}\theta\right) \theta^{2}w_{x}dx-\int \Lambda^{\alpha/2}\theta \left[\Lambda^{\alpha/2},w\right]\theta dx \\
& \leq C( \Vert \theta_{0} \Vert_{L^{\infty}}, \delta) \left(\Vert \sqrt{w}\theta \Vert_{L^{2}} \Vert \sqrt{w}\Lambda\theta \Vert_{L^{2}} + \Vert \sqrt{w} \theta \Vert^{2}_{L^{2}} \right) + \int \sqrt{w} \vert \Lambda^{\alpha/2}\theta \vert \ \frac{1}{\sqrt{w}}\left\vert \left[\Lambda^{\alpha/2},w\right]\theta \right \vert \ dx 
\end{split}
\end{equation*}

Then, using lemma \ref{lem24} one obtains, for $C( \Vert \theta_{0} \Vert_{L^{\infty}}, \delta)=C_1=(1+\delta)\Vert \theta_{0} \Vert_{L^{\infty}}$ for all $\eta_1>0$ and 
\begin{eqnarray*}
\frac{1}{2}\frac{d}{dt} \left\|\theta\right\|^{2}_{L^{2}(wdx)} + \left\|\Lambda^{\frac{\alpha}{2}}\theta\right\|^{2}_{L^{2}(wdx)} 
&\leq&  \left(\Vert\theta \Vert_{L^{2}(\w)} \Vert \Lambda\theta \Vert_{L^{2}_{w}} + \Vert \theta \Vert^{2}_{L^{2}({\w})} \right)C_1 + \Vert \theta \Vert_{L^{2}(\w)} \Vert \theta \Vert_{\dot H^{\alpha/2}(\w)}  \\
&\leq& C_2\Vert \theta \Vert^{2}_{ H^{1}(\w)} + \frac{1}{100}\Vert \theta \Vert^{2}_{\dot H^{\alpha/2}(\w)} 
\end{eqnarray*}

 Now, we study the contribution coming from the homogeneous Sobolev part  $\dot H^{1}(w)$. The transport part (corresponding to $\delta=0$) has been treated in \cite{Lazar2}, indeed, it is shown that
$$
\frac{1}{2} \partial_{t} \Vert \theta \Vert^{2}_{H^{1}(\w)} \leq C \Vert \theta \Vert^{2}_{H^{1}(\w)}.
$$
 Therefore, it just remains to estimate 
 $$
I_{\delta} \theta= -\delta \int w \Lambda \theta \Lambda (\theta \Lambda \theta) \ dx
 $$
 In order to take advantage of the dissipation, we start by rewritting $I_{\delta}$ as follows,
 \begin{eqnarray*}
 I_{\delta} \theta&=&-\delta \int \Lambda^{\mu}( w \Lambda \theta)  \Lambda^{1-\mu}(\theta \Lambda \theta) \ dx  \\
 &=& -\int \sqrt{w} \Lambda^{1-\mu}(\theta \Lambda \theta) \frac{1}{\sqrt{w}}[\Lambda^{\mu}, w] \Lambda \theta \ dx - \int \sqrt{w} \Lambda^{1+\mu} \theta \sqrt{w}   \Lambda^{1-\mu}(\theta \Lambda \theta) \ dx \\
 &=& I_1+I_2
  \end{eqnarray*}
 In order to estimate these two terms, we shall use the weighted Littlewood-Paley decomposition. We start with $I_1$, we have, by using  the commutator lemma \ref{lem24}  that
 $$
 \vert I_1 \vert \leq \Vert \theta \Lambda \theta \Vert_{\dot H^{1-\mu}(\w)} \Vert \theta \Vert_{\dot H^{1}(\w)}, 
 $$
 and then we estimate the right hand side by using the paraproduct formula \eqref{ppf} and we get
 $$
\Vert \theta \Lambda \theta \Vert_{\dot H^{1-\mu}(\w)} \leq \sum_{q \in \mathbb Z} 2^{q(1-\mu)} \Vert S_{q+1}\theta \Delta_{q} \mathcal{H}\theta_{x} \Vert_{L^{2}(\w)} +  \sum_{j \in \mathbb Z} 2^{j(1-\mu)} \Vert \Delta_{j}\theta S_{j} \mathcal{H}\theta_{x} \Vert_{L^{2}(\w)}.
$$
The first sum is controlled as follows
\begin{eqnarray*}
 \sum_{q\in \mathbb{Z}} 2^{q(1-\mu)} \Vert S_{q+1}\theta_{x}\ \Delta_{q} \mathcal{H}\theta_{x} \Vert_{L^{2}(\w)} &\leq&   \sum_{q \in \mathbb{Z}} 2^{q(1-\mu)} \Vert S_{q+1}\theta \Vert_{L^{\infty}}  \Vert \Delta_{q} \mathcal{H}\theta_{x} \Vert_{L^{2}(\w)} \\
 &\leq&  \Vert \theta\Vert_{L^{\infty}}   \sum_{q \in \mathbb{Z}} 2^{(1-\mu)q}   \Vert \Delta_{q} \mathcal \theta_{x} \Vert_{L^{2}(\w)}  \\
  &\leq&  \Vert \theta\Vert_{L^{\infty}}   \Vert \theta_{x} \Vert_{\dot H^{1-\mu}(\w)} \\
   &\leq&  \Vert \theta_{0}\Vert_{L^{\infty}}   \Vert \theta \Vert_{\dot H^{2-\mu}(\w)}
\end{eqnarray*}
where we used Berstein's inequality and the continuity of the Hilbert transform on $L^{2}(\w)$. \\

The other sum has to be treated in a different manner since we do not want to put $L^{\infty}$ in the term $S_{j} \mathcal{H} \theta_x$.   For instance, we can use H\"older's inequality and then Bernstein's inequality to recover some sufficiently nice Sobolev norm. 
\noindent We obtain, 
\begin{eqnarray*}
  \sum_{j\in \mathbb{Z}} 2^{j(1-\mu)} \Vert \Delta_{j}\theta S_{j} \mathcal{H}\theta_{x} \Vert_{L^{2}(w)} &\leq&   \sum_{j\in \mathbb{Z}} 2^{j(1-\mu)}  \Vert \Delta_{j}\theta \Vert_{L^{p}(w)} \Vert S_{j} \mathcal{H}\theta_{x} \Vert_{L^{q}(\w)} \\
  &\leq&  C   \sum_{j\in \mathbb{Z}} 2^{j(1-\mu)}   2^{j(\frac{1}{2}-\frac{1}{p})}  \Vert \Delta_{j}\theta \Vert_{L^{2}(w)} \Vert  \mathcal{H}\theta_{x} \Vert_{L^{q}(\w)}  \\
  &\leq&  C   \sum_{j\in \mathbb{Z}} 2^{j(1-\mu)}   2^{j(\frac{1}{2}-\frac{1}{p})}  \Vert \Delta_{j}\theta \Vert_{L^{2}(w)} 2^{j}\Vert  \mathcal{H}\theta \Vert_{L^{q}(\w)}  \\
    &\leq&  C  \Vert \theta \Vert_{L^{q}(\w)}  \sum_{j\in \mathbb{Z}} 2^{j(2-\mu+\frac{1}{q})} \Vert \Delta_{j}\theta \Vert_{L^{2}(\w)}  \\
  &\leq& C  \Vert \theta_{0} \Vert^{1-2/q}_{L^{\infty}} \Vert \theta \Vert^{2/q}_{L^{2}(w)} \Vert \theta \Vert_{\dot H^{2-\mu+\frac{1}{q}}(\w)}.
\end{eqnarray*}
%then using
%$$
%\Vert \theta \Vert_{\dot H^{5/4}(\w)} \leq \Vert \theta \Vert^{1/2}_{\dot H^{1}(\w)}\Vert \theta \Vert^{1/2}_{\dot H^{3/2}(\w)}
%$$
%$$
%\Vert \theta \Vert_{\dot H^{1/2}(\w)} \leq \Vert \theta \Vert^{1/2}_{L^{2}(\w)}\Vert \theta \Vert^{1/2}_{\dot H^{1}(\w)}
%$$
Hence,
\begin{eqnarray*}
 \vert I_1 \vert &\leq& \Vert \theta \Lambda \theta \Vert_{\dot H^{1-\mu}(\w)} \Vert \theta \Vert_{\dot H^{1}(\w)}  \\
 &\leq& \Vert \theta_{0}\Vert_{L^{\infty}}   \Vert \theta \Vert_{\dot H^{2-\mu}(\w)}\Vert \theta \Vert_{\dot H^{1}(\w)} 
+C \Vert \theta \Vert_{\dot H^{1}(\w)} \Vert \theta_{0} \Vert^{1-2/q}_{L^{\infty}} \Vert \theta \Vert^{2/q}_{L^{2}(w)} \Vert \theta \Vert_{\dot H^{2-\mu+\frac{1}{q}}(\w)}
 \end{eqnarray*}
 
 For $I_2$, following $I_1$, one gets
$$
\vert I_2 \vert  \leq \Vert \theta_{0}\Vert_{L^{\infty}}  \Vert \theta \Vert_{\dot H^{1+\mu}(\w)}  \Vert \theta \Vert_{\dot H^{2-\mu}(w)} +  \Vert \theta_{0} \Vert^{1-2/q}_{L^{\infty}} \Vert \theta \Vert_{\dot H^{1+\mu}(\w)} \Vert \theta \Vert^{2/q}_{L^{2}(\w)} \Vert \theta \Vert_{\dot H^{2-\mu+\frac{1}{q}}(\w)}.
$$
The idea is then to choose $q$ very small. Indeed, the associated Lebesgue space $L^{q}$ would be then a good substitute for $L^{\infty}$. It is worth recalling that this latter space is not well suited to estimate $\mathcal H \theta$ since it does not map $L^{\infty}$ to $L^{\infty}$ but $L^{\infty}$ to $BMO$. \\

Then, one has to choose $\mu$ close enough to $1/2$ because we control $H^{1+\frac{\alpha}{2}}(\w)$ with $\alpha \in (1,2)$. A good choice is to take for instance $\mu=\frac{1}{2}+\frac{1}{q}$, this choice allows us to interpolate and conclude the estimates. Indeed, for all $\delta \in (0,1)$
$$\Vert \theta \Vert_{\dot H^{2+\frac{1}{q}-\mu}(\w)} \leq 
\Vert \theta \Vert^{1-\delta}_{\dot H^{2+\frac{1}{q}-\mu-\delta}(\w)}
\Vert \theta \Vert^{\delta}_{\dot H^{3+\frac{1}{q}-\mu-\delta}(\w)},$$
Then, one observes that by choosing $\delta$  very close to 1 (essentially for the worst term, which is the second estimate), we get
$$
\Vert \theta \Vert^{1-\delta}_{\dot H^{2+\frac{1}{q}-\mu-\delta}(\w)}=\Vert \theta \Vert^{1-\delta}_{\dot H^{3/2-\delta}(\w)} \leq \Vert \theta \Vert^{1-\delta}_{ H^{1}(\w)}
$$
and if we choose $\delta=1-\epsilon$, with $\epsilon>0$ is chosen small enough so that $\frac{1}{2}+\epsilon \leq \frac{\alpha}{2}$, this is possible since $\alpha \in (1,2)$.
$$
\Vert \theta \Vert^{\delta}_{\dot H^{3+\frac{1}{q}-\mu-\delta}(\w)}=\Vert \theta \Vert^{\delta}_{\dot H^{3/2+ \epsilon}(\w)} \leq \Vert \theta \Vert^{\delta}_{H^{1+\frac{\alpha}{2}}(\w)}
$$
hence,
\begin{equation} \label{cte}
\Vert \theta \Vert_{\dot H^{2+\frac{1}{q}-\mu}(\w)} \leq  \Vert \theta \Vert^{1-\delta}_{ H^{1}}\Vert \theta \Vert^{\delta}_{H^{1+\frac{\alpha}{2}}(\w)}
\end{equation}
Therefore, since $\frac{1}{2}+\frac{1}{q} \leq \frac{\alpha}{2}$, one obtains

\begin{eqnarray} \label{sim}
 \vert I_1 \vert  &\lesssim&   \Vert \theta \Vert_{H^{\frac{1}{2}+\frac{\alpha}{2}}(\w)}\Vert \theta \Vert_{\dot H^{1}(\w)} 
+  \Vert \theta \Vert^{2/q}_{L^{2}(\w)}  \Vert \theta \Vert^{2-\delta}_{ H^{1}(\w)}\Vert \theta \Vert^{\delta}_{H^{1+\frac{\alpha}{2}}(\w)}
 \end{eqnarray}
Then, using Young's inequality with $p_1=\frac{2-\nu}{2-\delta}$ and its conjugate $p_2=\frac{2-\nu}{\delta-\nu}$, with $0<\nu<\delta$ where $\nu$ will be chosen later. One gets, for all $\epsilon_1>0$
\begin{eqnarray*}
 \vert I_1 \vert  &\lesssim&  \frac{\epsilon_1}{2} \Vert \theta \Vert^{2}_{H^{\frac{1}{2}+\frac{\alpha}{2}}(\w)} + \frac{1}{2\epsilon_1}\Vert \theta \Vert^{2-\nu}_{\dot H^{1}(\w)}  
+   \frac{1}{2\epsilon_1}\Vert \theta \Vert^{\frac{2}{q} \frac{2-\nu}{2-\delta}}_{L^{2}(w)}  \Vert \theta \Vert^{2-\nu}_{ H^{1}(\w)}
+  \frac{\epsilon_1}{2}\Vert \theta \Vert^{\frac{\delta}{2-\delta}(2-\nu)}_{H^{1+\frac{\alpha}{2}}(\w)} 
 \end{eqnarray*}
 And once again, we use Young's inequality with $p=\frac{2}{2-\nu}$ hence $q=\frac{\nu}{2}$
\begin{eqnarray*}
 \vert I_1 \vert  &\lesssim&  \frac{\epsilon_1}{2} \Vert \theta \Vert^{2}_{H^{\frac{1}{2}+\frac{\alpha}{2}}(\w)} + \frac{1}{2\epsilon_1}\Vert \theta \Vert^{2}_{\dot H^{1}(\w)}  
+   \frac{C_1(\delta,\nu)}{2\epsilon_1}\Vert \theta \Vert^{\frac{4}{q} \frac{1}{2-\delta} }_{L^{2}(\w)}  + \frac{C_2(\delta,\nu)}{2\epsilon_1} \Vert \theta \Vert^{(2-\nu)\frac{\nu}{2}}_{ H^{1}(\w)}
+  \frac{\epsilon_1}{2}\Vert \theta \Vert^{\frac{\delta}{2-\delta}(2-\nu)}_{H^{1+\frac{\alpha}{2}}(\w)} 
 \end{eqnarray*}
 Since $q$ is chosen large, then $\frac{4}{q} \frac{1}{2-\delta}\leq2$, we also have that $(2-\nu)\frac{\nu}{2}\leq2$ (this is a consequence of $(\nu-2)^2\geq0$). We can choose any value of $\nu \in (0,\delta)$
so that $\frac{\delta}{2-\delta}(2-\nu)\leq 2$. For instance, if $\nu=1/2$ then we obviously have $\frac{\delta}{2-\delta}\leq 4/3$ for all $\delta \in (0,1)$ (importantly, there is no restriction on $\delta$, the previous estimates hold for all $\delta \in(0,1)$, we shall use this fact to control the second term, see \ref{sim2} below).

One finally obtains
\begin{eqnarray} \label{I1}
 \vert I_1 \vert  \lesssim  {\epsilon_1} \Vert \theta \Vert^{2}_{H^{1+\frac{\alpha}{2}}(\w)} + \frac{2}{\epsilon_1}\Vert \theta \Vert^{2}_{ H^{1}(\w)}  
 \end{eqnarray}
 For $I_2$, we have seen that
 $$
\vert I_2 \vert  \leq \Vert \theta_{0}\Vert_{L^{\infty}}  \Vert \theta \Vert_{\dot H^{1+\mu}(\w)}  \Vert \theta \Vert_{\dot H^{2-\mu}(\w)} +  \Vert \theta_{0} \Vert^{1-2/q}_{L^{\infty}} \Vert \theta \Vert_{\dot H^{1+\mu}(\w)} \Vert \theta \Vert^{2/q}_{L^{2}(\w)} \Vert \theta \Vert_{\dot H^{2-\mu+\frac{1}{q}}(\w)}.
$$
Since $\mu=\frac{1}{2}+\frac{1}{q}$ and $q$ is such that $\frac{1}{2}+\frac{1}{q} \leq \frac{\alpha}{2}$ then by using \ref{cte}, the latter inequality becomes
\begin{eqnarray*}
\vert I_2 \vert  &\leq& \Vert \theta_{0}\Vert_{L^{\infty}}  \Vert \theta \Vert_{\dot H^{\frac{3}{2}+\frac{1}{q}}(\w)}  \Vert \theta \Vert_{\dot H^{\frac{3}{2}-\frac{1}{q}}(\w)} +  \Vert \theta_{0} \Vert^{1-2/q}_{L^{\infty}} \Vert \theta \Vert_{\dot H^{\frac{3}{2}+\frac{1}{q}}(w)} \Vert \theta \Vert^{2/q}_{L^{2}(\w)} \Vert \theta \Vert^{1-\delta}_{ H^{1}(\w)}\Vert \theta \Vert^{\delta}_{H^{1+\frac{\alpha}{2}}(\w)} \\
&\leq&\Vert \theta_{0}\Vert_{L^{\infty}}  \Vert \theta \Vert_{\dot H^{1+\frac{\alpha}{2}}(\w)}  \Vert \theta \Vert_{\dot H^{\frac{3}{2}-\frac{1}{q}}(w)} +  \Vert \theta_{0} \Vert^{1-2/q}_{L^{\infty}} \Vert \theta \Vert_{\dot H^{\frac{3}{2}+\frac{1}{q}}(\w)} \Vert \theta \Vert^{2/q}_{L^{2}(\w)} \Vert \theta \Vert^{1-\delta}_{ H^{1}(\w)}\Vert \theta \Vert^{\delta}_{H^{1+\frac{\alpha}{2}}(\w)} \\
&\lesssim& \Vert \theta \Vert^{2-\frac{1}{q}}_{ H^{1+\frac{\alpha}{2}}(\w)}   \Vert \theta \Vert^{\frac{1}{q}}_{ H^{1}(\w)} +   
\Vert \theta \Vert^{1+\delta}_{H^{1+\frac{\alpha}{2}}(\w)} \Vert \theta \Vert^{2/q}_{L^{2}(\w)} \Vert \theta \Vert^{1-\delta}_{ H^{1}(\w)}
 \\
\end{eqnarray*}
where, in the last step, we used the following  interpolation inequality 
$$
\Vert \theta \Vert_{ H^{\frac{3}{2}-\frac{1}{q}}(\w)} \leq \Vert \theta \Vert^{1-\frac{1}{q}}_{ H^{\frac{3}{2}}(\w)}  \Vert \theta \Vert^{\frac{1}{q}}_{ H^{\frac{1}{2}}(\w)},
$$
Then, we use Young's inequality with $p_3=2q$ and its conjugate $p_4=\frac{2q}{2q-1}$ in the first product 
\begin{eqnarray*}
 \Vert \theta \Vert^{2-\frac{1}{q}}_{ H^{1+\frac{\alpha}{2}}(\w)}   
 \Vert \theta \Vert^{\frac{1}{q}}_{ H^{1}(\w)}  &\lesssim& \epsilon_{1}\Vert \theta \Vert^{2}_{ H^{1+\frac{\alpha}{2}}(\w)} + \frac{1}{\epsilon_{1}}\Vert \theta \Vert^{2}_{ H^{1}(\w)}
 \end{eqnarray*}
 it remains to estimate
 $$
 L=\Vert \theta \Vert^{1+\delta}_{H^{1+\frac{\alpha}{2}}(\w)} \Vert \theta \Vert^{2/q}_{L^{2}(\w)} \Vert \theta \Vert^{1-\delta}_{ H^{1}(\w)}
 $$
 but $\delta=1-\epsilon,$ so that 
 \begin{equation} \label{sim2}
 L=\Vert \theta \Vert^{2-\epsilon}_{H^{1+\frac{\alpha}{2}}(\w)} \Vert \theta \Vert^{2/q}_{L^{2}(\w)} \Vert \theta \Vert^{\epsilon}_{ H^{1}(\w)}
\end{equation}
 but, we have already controlled the same term with $\delta>0$ in stead of $\epsilon>0$ (see the second term in the right hand side of \ref{sim}), therefore following the same steps, one arrives at the same conclusion at in $I_1$ that is to say, for $i=1,2$, one has
 
 \begin{eqnarray} \label{I2}
 \vert I_i \vert  \lesssim  {\epsilon_1} \Vert \theta \Vert^{2}_{H^{1+\frac{\alpha}{2}}(\w)} + \frac{2}{\epsilon_1}\Vert \theta \Vert^{2}_{ H^{1}(\w)}  
 \end{eqnarray}
 Hence, to conclude it suffices to choose $\epsilon_1$ sufficiently small. 
 
  \qed
 
 \section{Proof of Theorem \ref{super}}

 The control  of the $L^{2}(w)$ norm is straightforward, indeed, using the    lemma \ref{lem24} along with the weighted Sobolev embedding $\dot H^{1/6}({w}) \hookrightarrow L^{3}({w})$, one obtains
  \begin{eqnarray*} 
\partial_{t} \int \frac{\theta^{2}}{2}w \ dx + \int \vert \Lambda^{\alpha/2} \theta \vert^{2} w \ dx &=& - \int \theta^{2} \Lambda \theta w \ dx \ + \int \theta^{2} \mathcal{H}\theta w_x \ dx +\delta \int w \theta^{2} \Lambda \theta \ dx \\
&& \ -  \nonumber \int \Lambda^{\alpha/2} \theta [\Lambda^{\alpha/2},w] \theta \ dx \\
&\lesssim& \Vert \theta \Vert^{2}_{\dot H^{1/6}({\w})}\Vert \theta \Vert_{\dot H^{7/6}({\w})}+\Vert \theta \Vert^{3}_{\dot H^{1/6}({\w})} + \Vert \theta \Vert^{2}_{H^{\alpha/2}(\w)} \lesssim \Vert \theta \Vert^{3}_{H^{2}(\w)}
\end{eqnarray*}

The control of the homogeneous part $\dot H^{2}(\w)$ is  done as follows.
\begin{eqnarray*}
\frac{1}{2}\partial_{t} \Vert \theta \Vert^{2}_{H^{2}(\w)} + \int  \vert \Lambda^{2+\frac{\alpha}{2}} \theta \vert^{2} \ w \ dx &=&  -\frac{1}{2}\int w_x (\theta_{xx})^{2}  \mathcal{H} \theta  -\frac{3}{2}\int w(\theta_{xx})^{2} \Lambda \theta  + \int w \theta_{xx}  \theta_{x} \Lambda \theta_{x} \\
&&- \int \Lambda^{2+\frac{\alpha}{2}} \theta [\Lambda^{\alpha/2},w] \theta_{xx} \\
&-& \delta \int w(\theta_{xx})^2 \Lambda \theta \ dx -2\delta \int w\theta_{xx}\theta_{x} \Lambda \theta_{x} \ dx
-\delta \int w\theta\theta_{xx} \Lambda \theta_{xx} \ dx \\ 
&=&\sum_{j=1}^{7} I_j
\end{eqnarray*}
We already know from \cite{Lazar2} that
\begin{equation*} 
\left\vert  \sum_{j=1}^{4} I_j \right \vert \lesssim  \frac{1}{\epsilon_1^{4} } \Vert \theta \Vert^{4}_{ H^{2}({\w})} + {\epsilon_1}\Vert \theta \Vert^{\frac{16}{3}}_{H^{2}(\w) }+ {\epsilon_1}\Vert \theta \Vert^{2}_{\dot H^{2+\frac{\alpha}{2}}({\w})}
\end{equation*}
where we crucially used the first point of lemma \ref{lem24} to estimate $I_4$. Futhermore, one observes that $I_5$ and $I_6$ are the same as respectively $I_2$ and $I_3$, hence
\begin{equation} \label{ineg}
\left\vert  \sum_{j=1}^{6} I_j \right \vert \lesssim  \frac{1}{\epsilon_1^{4} } \Vert \theta \Vert^{4}_{ H^{2}({\w})} + {\epsilon_1}\Vert \theta \Vert^{\frac{16}{3}}_{H^{2}(\w) }+ {\epsilon_1}\Vert \theta \Vert^{2}_{\dot H^{2+\frac{\alpha}{2}}({\w})}
\end{equation}
 \\

\noindent The most singular term is $I_7$. Since $\delta\geq0$ and $\theta\geq0$ we may use the pointwise inequality \ref{pointwise lemma} and then by integrating by parts we get  
 \begin{eqnarray*}
 I_7&=& -\delta \int w\theta\theta_{xx} \Lambda \theta_{xx} \ dx \\
 &\leq&-\delta \int w\theta \Lambda ((\theta_{xx})^2)  \ dx \\
 &=&-\delta \int \Lambda (w\theta) \ (\theta_{xx})^2  \ dx  \\
 \end{eqnarray*}
 At this step, it is important to note that one cannot take advantage of the $L^{\infty}$ maximum principle, indeed, this would imply to estimate the remaning term in $L^{1}(w)$ and it is well-known that the Hilbert transform is not continuous on that space. The idea here is to write the term $\Lambda(w\theta)$ as a controlled commutator (via lemma \ref{lem24}) plus another term which is easy to control since the weight will be outside the differential term. More precisely, we write
  \begin{eqnarray*}
 I_7&=&-\delta \int (\theta_{xx})^2 [\Lambda, w] \theta \ dx - \delta \int w(\theta_{xx})^2 \Lambda \theta \ dx  \\
% &=&-\delta \int {w}(\theta_{xx})^2 \ \frac{1}{{{w}}} [\Lambda,{w}]\theta \ dx \ + \ \delta \int {w}(\theta_{xx})^2  \Lambda \theta \ dx
% \\
 &=&I_{7,1} + I_{7,2}
\end{eqnarray*}
To estimate $I_{7,1}$ we use H\"older's inequality with $\frac{\alpha}{2} + \frac{1-\alpha}{2}=\frac{1}{2}$, and the two following weighted Sobolev embeddings $\dot H^{\alpha/2} (\w) \hookrightarrow L^{\frac{2}{1-\alpha}}(\w)$  and  $\dot H^{\frac{1-\alpha}{2}} (\w) \hookrightarrow L^{\frac{2}{\alpha}}(\w)$, and the last point of lemma \ref{lem24}, one gets
 \begin{eqnarray*}
\vert I_{7,1} \vert&\leq& \Vert \theta_{xx} \Vert_{L^{2}(\w)} \Vert \theta_{xx}\Vert_{L^{\frac{2}{1-\alpha}}(\w)} \left\Vert  \frac{1}{w}[\Lambda, w] \theta \right\Vert_{L^{\frac{2}{\alpha}}(\w)} \\
&\leq& \Vert \theta \Vert_{\dot H^2 (\w)} \Vert \theta_{xx}\Vert_{\dot H^{\alpha/2}(\w)} \Vert \theta \Vert_{L^{\frac{2}{\alpha}}(\w)} \\
&\leq& \Vert \theta \Vert_{\dot H^2 (\w)}\Vert \theta \Vert_{\dot H^{2+\alpha/2}(\w)}\Vert \theta
 \Vert_{\dot H^{\frac{1-\alpha}{2}} (\w)}   \\
&\leq& \Vert \theta \Vert_{\dot H^2 (\w)}\Vert \theta \Vert_{\dot H^{2+\alpha/2}(\w)}\Vert \theta
 \Vert_{ H^{\frac{1}{2}} (\w)} \\
&\leq& \frac{1}{\epsilon_1} \Vert \theta \Vert^{4}_{ H^{2}(\w)} + \epsilon_1 \Vert \theta \Vert^{2}_{\dot H^{2+\frac{\alpha}{2}} (\w)}
  \end{eqnarray*}
  For $I_{7,2}$ we have, for all $\epsilon_1>0$,
  \begin{eqnarray*}
\vert I_{7,2} \vert &\leq& \Vert w^{\frac{1}{2}}\theta_{xx} \Vert_{L^{2}} \Vert w^{\frac{1-\alpha}{2}} \theta_{xx} \Vert_{L^{\frac{2}{1-\alpha}}} \left\Vert w^{\frac{\alpha}{2}}\Lambda\theta \right \Vert_{L^{2/\alpha}} \\
&\leq& \Vert \theta \Vert_{\dot H^2 (\w)}\Vert \theta \Vert_{\dot H^{2+\frac{\alpha}{2}} (\w)}  \Vert \theta \Vert_{\dot H^{\frac{3-\alpha}{2}} (\w)} \\
&\leq& \frac{1}{\epsilon_1} \Vert \theta \Vert^{4}_{ H^{2}(\w)} + \epsilon_1 \Vert \theta \Vert^{2}_{\dot H^{2+\frac{\alpha}{2}} (\w)}
  \end{eqnarray*}
  Therefore, 
  \begin{equation} \label{esti}
   \vert I_7 \vert \leq \frac{2}{\epsilon_1} \Vert \theta \Vert^{4}_{ H^{2}(\w)} + 2\epsilon_1 \Vert \theta \Vert^{2}_{\dot H^{2+\frac{\alpha}{2} }(\w)}
  \end{equation}
  Combining \ref{ineg} and \ref{esti} we conclude by choosing $\epsilon_1>0$ 
 \begin{equation} \label{ineg}
\left\vert  \sum_{j=1}^{4} I_j \right \vert \lesssim  \frac{1}{\epsilon_1^{4} } \Vert \theta \Vert^{4}_{ H^{2}({\w})} + {\epsilon_1}\Vert \theta \Vert^{\frac{16}{3}}_{H^{2}(\w) }+ {\epsilon_1}\Vert \theta \Vert^{2}_{\dot H^{2+\frac{\alpha}{2}}({\w})}
\end{equation}
hence, by choosing $\epsilon_1>0$ suficiently small, one finds  
  \begin{eqnarray*}
\frac{1}{2}\partial_{t} \Vert \theta \Vert^{2}_{H^{2}{\w}} + C\int  \vert \Lambda^{2+\frac{\alpha}{2}} \theta \vert^{2} \ w \ dx &\lesssim& \Vert\theta \Vert^{4}_{ H^{2}({\w})} + \Vert \theta \Vert^{\frac{16}{3}}_{H^{2}(\w) } 
  \end{eqnarray*}
  By using Gr\"onwall's inequality, and the classical truncation and compactness arguments (see \cite{Lazar2} for instance), one gets the theorem. 
  \qed

%%%%%%%%%%%%%%%%%%%%%%%%%
\section{Proof of Theorem \ref{GW Wiener}}
%%%%%%%%%%%%%%%%%%%%%%%%%

%%%%%%%%%%%%%%%%%
\subsection{A priori estimates}
%%%%%%%%%%%%%%%%%

Taking the Fourier transform of \eqref{transport equation 1}, we have that
\begin{align*}
\pat|\hat{\theta}(\xi)| &=\frac{\bar{\hat{\theta}}(\xi)\pat\hat{\theta}(\xi)+\hat{\theta}(\xi)\pat\bar{\hat{\theta}}(\xi)}{2|\hat{\theta}(\xi)|}=\frac{\text{Re}\left(\bar{\hat{\theta}}(\xi)\pat\hat{\theta}(\xi)\right)}{|\hat{\theta}(\xi)|}\\
&=-\text{Re}\left[\int_{\RR}\frac{-i\zeta}{|\zeta|}\hat{\theta}(\zeta)i(\xi-\zeta))\hat{\theta}(\xi-\zeta)-\delta\hat{\theta}(\zeta)|\xi-\zeta|\hat{\theta}(\xi-\zeta)d\zeta\frac{\bar{\hat{\theta}}(\xi)}{|\hat{\theta}(\xi)|}\right]\frac{1}{\sqrt{2\pi}}-|\xi||\hat{\theta}|.
\end{align*}
Consequently,
\begin{align*}
\frac{d}{dt}\|\theta\|_{A^0} &\leq (1+|\delta|) \int_{\RR}\int_{\RR}|\hat{\theta}(\zeta)||\xi-\zeta||\hat{\theta}(\xi-\zeta)|d\zeta d\xi \frac{1}{\sqrt{2\pi}}-\|\theta\|_{\dot{A}^1}\\&\leq (1+|\delta|) \int_{\RR}\int_{\RR}|\hat{\theta}(\zeta)||\xi-\zeta||\hat{\theta}(\xi-\zeta)| d\xi d\zeta \frac{1}{\sqrt{2\pi}}-\|\theta\|_{\dot{A}^1}\leq \left(\frac{(1+|\delta|)\|\theta\|_{A^0}}{\sqrt{2\pi}}-1\right)\|\theta\|_{\dot{A}^1}.
\end{align*}
Thus, if $\theta_{0}$ satisfies the condition (\ref{initial data in Wiener}), we have
\begin{equation}\label{decay1}
\|\theta(t)\|_{A^0}+\left(1-\frac{(1+|\delta|)\|\theta_0\|_{A^0}}{\sqrt{2\pi}}\right)\int_0^\infty \|\theta(s)\|_{\dot{A}^1}ds\leq \|\theta_0\|_{A^0}.
\end{equation}
Similarly, 
\begin{align*}
\pat|\xi \hat{\theta}| &=-\text{Re}\bigg{[}\int_{\RR}\frac{-i\zeta}{|\zeta|}\hat{\theta}(\zeta)(i(\xi-\zeta))^2\hat{\theta}(\xi-\zeta)+|\zeta|\hat{\theta}(\zeta)i(\xi-\zeta))\hat{\theta}(\xi-\zeta)d\zeta\\
&\quad-\delta\int_{\RR}i\zeta\hat{\theta}(\zeta)|\xi-\zeta|\hat{\theta}(\xi-\zeta)+\hat{\theta}(\zeta)i(\xi-\zeta)|\xi-\zeta|\hat{\theta}(\xi-\zeta)d\zeta\frac{\overline{i\xi\hat{\theta}}(\xi)}{|\xi||\hat{\theta}(\xi)|}\bigg{]}\frac{1}{\sqrt{2\pi}}-|\xi|^2|\hat{\theta}|.
\end{align*}
Thus, using \eqref{eq:5d}, we have that
\begin{align*}
\frac{d}{dt}\|\theta\|_{\dot{A}^1} &\leq \left(\|\theta\|_{A^0}\|\theta\|_{\dot{A}^2}+\|\theta\|_{A^1}^2\right)\frac{1+|\delta|}{\sqrt{2\pi}}-\|\theta\|_{\dot{A}^2}\leq \left(\frac{2(1+|\delta|)\|\theta\|_{A^0}}{\sqrt{2\pi}} -1\right)\|\theta\|_{\dot{A}^2}.
 \end{align*}
As a consequence, we obtain that, if  $\theta_{0}$ satisfies the condition (\ref{initial data in Wiener}), we also have  
\begin{equation}\label{decay2}
\|\theta(t)\|_{\dot{A}^1}+\left(1-\frac{2(1+|\delta|)\|\theta_{0}\|_{A^0}}{\sqrt{2\pi}}\right)\int_0^\infty \|\theta(s)\|_{\dot{A}^2}ds\leq \|\theta_0\|_{\dot{A}^1}.
\end{equation}
By (\ref{decay1}) and (\ref{decay2}), we conclude that 
\eqn \label{decay3}
\|\theta(t)\|_{A^1}+\left(1-\frac{2(1+|\delta|)\|\theta_{0}\|_{A^0}}{\sqrt{2\pi}}\right)\int_0^t \|\theta_{x}(s)\|_{A^1}ds\leq \|\theta_0\|_{A^1}
\een
for all $t\ge 0$.

%%%%%%%%%%%%%%%%%%%%%%%%%%%%
\subsection{Approximation and passing to the limit}
%%%%%%%%%%%%%%%%%%%%%%%%%%%%
Define $e^{\epsilon\pax^2}$ the heat semigroup, \emph{i.e.}
\[
\widehat{e^{\epsilon\pax^2} f(x)}=e^{-\epsilon \xi^2}\hat{f}(\xi),
\]
and $g_\epsilon(x)=e^{-\epsilon x^2}$. Note that
\[
\hat{g}_\epsilon(\xi)=\frac{1}{\sqrt{2\epsilon}} e^{-\frac{\xi^2}{4\epsilon}}, \quad \left\|\hat{g}_\epsilon\right\|_{L^1}=\sqrt{2\pi}.
\]
Given $\theta_0(x)\in A^1$, we consider $\theta^\epsilon_0(x)=g_\epsilon(x)e^{\epsilon\pax^2}\theta_0(x)$. As $\theta_0(x)$ is a bounded function, we have that $\theta_0^\epsilon$ is infinitely smooth and has finite total mass:
\[
\left\|\theta_0^\epsilon\right\|_{L^1}\leq \|g_\epsilon\|_{L^1} \left\|e^{\epsilon\pax^2}\theta_0 \right\|_{L^\infty}\leq \sqrt{\frac{\pi}{\epsilon}}\|\theta_0\|_{L^\infty}.
\]
Furthermore, using Young's inequality and the definition of $g_\epsilon$,
\[
\left\|\theta_0^\epsilon\right\|_{A^0}=\frac{1}{\sqrt{2\pi}} \left\|\hat{g}_\epsilon*\left(e^{-\epsilon \xi^2}\hat{\theta_0}\right)\right\|_{L^1}\leq \left\|e^{-\epsilon \xi^2}\hat{\theta}_0 \right\|_{L^1}\leq \|\theta_0\|_{A^0}.
\]
Similarly,
\[
\|\theta_0^\epsilon\|_{\dot{A}^1}=\|\pax \theta_0^\epsilon\|_{A^0}\leq \|\theta_0\|_{\dot{A}^1}+\left\|\pax g_\epsilon e^{\epsilon\pax^2}\theta_0\right\|_{A^0}=\|\theta_0\|_{\dot{A}^1}+c\epsilon\|\theta_0\|_{A^0}.
\]

Now we define the approximated problems 
\eqn \label{reg equation}
\theta^\epsilon_t+\left(\mathcal{H}\theta^\epsilon\right)\theta^\epsilon_x + \delta\theta^\epsilon \Lambda\theta^\epsilon+\Lambda\theta^\epsilon=\epsilon\pax^2\theta^\epsilon,
\een
with finite energy approximated initial data $\theta_0^\epsilon$. These problems have unique smooth solutions denoted by $\theta^\epsilon$. Moreover, $(\theta^{\epsilon})$ satisfies a uniform bound 
\[
\|\theta^{\epsilon}(t)\|_{A^1}+\left(1-\frac{2(1+|\delta|)\|\theta_{0}\|_{A^0}}{\sqrt{2\pi}}\right)\int_0^t \|\theta^{\epsilon}_{x}(s)\|_{A^1}ds +\epsilon \int_0^t \|\theta^{\epsilon}_{x}(s)\|_{A^2}ds\leq \|\theta_0\|_{A^1}.
\]
uniformly in $\epsilon$. Thus, the a priori estimates lead to the following uniform-in-$\epsilon$ bounds
\begin{equation}\label{uniform bound Wiener}
\begin{split}
& \sup_{t\in[0,\infty)}\|\theta^\epsilon(t)\|_{C^0(\RR)}\leq \sup_{t\in[0,\infty)}\|\theta^\epsilon(t)\|_{A^0}\leq \|\theta_0\|_{A^0}<\frac{\sqrt{2\pi}}{1+|\delta|},\\
& \sup_{t\in[0,\infty)}\|\theta^\epsilon(t)\|_{\dot{C}^1(\RR)}\leq \sup_{t\in[0,\infty)}\|\theta^\epsilon(t)\|_{\dot{A}^1}<\|\theta_0\|_{\dot{A}^1}+c\epsilon\|\theta_0\|_{A^0},\\
& \sup_{t\in[0,\infty)}\|\Lambda\theta^\epsilon(t)\|_{C^0(\RR)}<\|\theta_0\|_{\dot{A}^1}+c\epsilon\|\theta_0\|_{A^0}.
\end{split}
\end{equation}
Moreover, from the equation (\ref{reg equation}) we also obtain uniform bounds
\begin{equation}
\begin{split}
&\sup_{t\in[0,\infty)}\|\pat\theta^\epsilon(t)\|_{C^0(\RR)}+\sup_{t\in[0,\infty)}\|\pat H\theta^\epsilon(t)\|_{C^0(\RR)}\leq 2\sup_{t\in[0,\infty)}\|\pat\theta^\epsilon(t)\|_{A^0(\RR)}\leq F_1(\|\theta_0\|_{A^1},\delta),\\
& \|H\theta^\epsilon\|_{C^1([0,\infty)\times \RR)}+\|\theta^\epsilon\|_{C^1([0,\infty)\times \RR)}\leq F_2(\|\theta_0\|_{A^1},\delta)
\end{split}
\end{equation}
where $F_{1}$ and $F_{2}$ only depend on the quantities in the right-hand side of (\ref{uniform bound Wiener}). Due to Banach-Alaoglu Theorem, there exists a subsequence (denoted by $\epsilon$) and a limit function $\theta\in W^{1,\infty}([0,\infty)\times\RR)$ such that 
\[
\theta^\epsilon\stackrel{\star}{\rightharpoonup} \theta \quad \text{in} \quad L^{\infty}\left(0,T; W^{1,\infty}\right), \quad \Lambda\theta^\epsilon\stackrel{\star}{\rightharpoonup} \Lambda \theta \quad \text{in} \quad L^{\infty}\left(0,T; L^{\infty}\right)
\]
for all $T>0$. Using Lemma \ref{IVlemapaco}, we have the following strong convergence
\[
\lim_{\epsilon\rightarrow0}\|\theta^{\epsilon}-\theta\|_{L^\infty(K)}+\|H\theta^{\epsilon}-H\theta\|_{L^\infty(K)}=0
\]
for $K=[0,T]\times[-R,R]$ for any $R>0$. Since $L^{\infty}(K)\subset L^{2}(K)$, we can pass to the limit to the weak formulation
\[
\int^{T}_{0}\int \left[\theta^{\epsilon} \psi_{t} + \left(\mathcal{H}\theta^{\epsilon}\right)\theta^{\epsilon}\psi_{x}+(1-\delta)\Lambda\theta^{\epsilon}\theta^{\epsilon} \psi -\theta^{\epsilon}\Lambda\psi +\epsilon \theta^{\epsilon}\psi_{xx}\right]  dxdt =\int\theta^{\epsilon}_{0}(x)\psi(0,x)dx
\]
to obtain a weak solution $\theta$ satisfying 
\[
\theta\in\mathcal{W}_T, \quad \sup_{t\in[0,T]}\|\theta(t)\|_{A^1}+\left(1-\frac{\sqrt{2} \|\theta_{0}\|_{A^{0}}}{\sqrt{\pi}}\right)\int^{T}_{0}\|\theta_{x}(t)\|_{A^1} dt\leq \|\theta_0\|_{A^1}
\]
for all $T>0$.

%%%%%%%%%%%%%%%
\subsection{Uniqueness}
%%%%%%%%%%%%%%%
To show the uniqueness of a weak solution, we consider the equation of $\theta=\theta_{1}-\theta_{2}$ given by 
\eqn \label{difference Fourier}
\theta_{t}+\Lambda \theta=-\left(\mathcal{H}\theta\right) \theta_{1x}- \left(\mathcal{H}\theta_{2}\right) \theta_{x}-\delta \theta \Lambda \theta_{1}-\delta \theta_{2}\Lambda \theta,  \quad \theta(0,x)=0.  
\een
Taking the Fourier transform of \eqref{difference Fourier} and multiply by $\frac{\hat{\theta}}{|\hat{\theta}|}$, we have
\begin{equation}
\begin{split}
\frac{d}{dt}\left\|\theta\right\|_{A^{0}}+\|\theta\|_{\dot{A}^{1}}\leq \frac{1+|\delta|}{\sqrt{2\pi}}\left\|\theta_{1}\right\|_{\dot{A}^{1}} \|\theta\|_{A^{0}}+\frac{1+|\delta|}{\sqrt{2\pi}}\left\|\theta_{2}\right\|_{A^{0}} \|\theta\|_{\dot{A}^{1}}.
\end{split}
\end{equation}
Since
\[
\frac{(1+|\delta|)\left\|\theta_{2}\right\|_{A^{0}}}{\sqrt{2\pi}}<1, \quad \left\|\theta_{1}(t)\right\|_{\dot{A}^{1}} \in L^{1}([0,\infty)),
\]
we have $\theta(t,x)=0$ in $A^{0}$ for all time. This implies that a weak solution is unique.

%%%%%%%%%%%%%%%%%%%%%
\section{Proof of Theorem \ref{GW 2}}
%%%%%%%%%%%%%%%%%%%%%

%%%%%%%%%%%%%%%%
\subsection{A priori estimate}
%%%%%%%%%%%%%%%%
We consider the equation
\eqn \label{equation 2}
\theta_{t}+u\theta_{x}+\Lambda^{\gamma}\theta=0, \quad u=(1-\partial_{xx} )^{-\alpha} \theta
\een
Since (\ref{equation 2}) satisfies the maximum principle, we have 
\[
\|\theta(t)\|_{L^\infty}\leq \|\theta_0\|_{L^\infty}.
\]
We also obtain that 
\begin{equation*}
\begin{split}
\frac{1}{2}\frac{d}{dt}\left\|\theta\right\|^{2}_{L^{2}}+\left\|\Lambda^{\frac{\gamma}{2}}\theta\right\|^{2}_{L^{2}}&=-\int u\theta_{x}\theta dx\leq \|\theta\|_{L^{\infty}}\|u_{x}\|_{L^{2}}\|\theta\|_{L^{2}}\\
&\leq C\left(\|\theta_{0}\|_{L^{\infty}}+ \|\theta_{0}\|^{2}_{L^{\infty}}\right)\|\theta\|^{2}_{L^{2}} +\frac{1}{2}\left\|\Lambda^{\frac{\gamma}{2}}\theta\right\|^{2}_{L^{2}}
\end{split}
\end{equation*}
where we use the condition $\alpha= \frac{1}{2}-\frac{\gamma}{4}$ to bound $u_{x}$ as  
\[
\|u_{x}\|_{L^{2}}\leq C\left(\|\theta\|_{L^{2}} + \left\|\Lambda^{\frac{\gamma}{2}}\theta\right\|_{L^{2}}\right).
\]
Therefore, we obtain that
\[
\left\|\theta(t)\right\|^{2}_{L^{2}}+\int^{t}_{0}\left\|\Lambda^{\frac{\gamma}{2}}\theta(s)\right\|^{2}_{L^{2}}ds \leq \|\theta_{0}\|^{2}_{L^{2}}e^{C\left(1+ \|\theta_{0}\|^{2}_{L^{\infty}}\right)t}.
\]

%%%%%%%%%%%%%%%%%%%%%%%%%%
\subsection{Approximation and passing to limit}
%%%%%%%%%%%%%%%%%%%%%%%%%%
We consider the following equation with regularized initial data: 
\eqn \label{regularized equation 2}
\theta^{\epsilon}_t+u^{\epsilon}\theta^{\epsilon}_x+\Lambda^{\gamma}\theta^{\epsilon}=0,\quad  \theta^{\epsilon}_{0}=\rho_{\epsilon}\ast \theta_{0}.
\een
Then, there exists a global-in-time smooth solution $\theta^{\epsilon}$. Moreover, $\theta^{\epsilon}$ satisfies that
\begin{equation} 
\begin{split}
&\left\| \theta^{\epsilon}(t) \right\|_{L^{\infty}}+ \left\| \theta^{\epsilon}(t) \right\|^{2}_{L^{2}} +\int^{t}_{0}\left\|\Lambda^{\frac{\gamma}{2}}\theta^{\epsilon}(s)\right\|^{2}_{L^{2}}ds +\epsilon\left\|\nabla \theta^{\epsilon}\right\|^{2}_{L^{2}} \leq \left\| \theta_{0}\right\|_{L^{\infty}}+ \|\theta_{0}\|^{2}_{L^{2}}e^{C\left(1+ \|\theta_{0}\|^{2}_{L^{\infty}}\right)t}.
\end{split}
\end{equation}
Therefore, $(\theta_{\epsilon})$ is bounded in $\mathcal{C}_{T}$ uniformly in $\epsilon>0$. This implies the uniform bounds 
\[
u^{\epsilon}\in L^{4}\left(0,T; L^{4}\right), \quad \theta^{\epsilon}\in L^{4}\left(0,T; L^{4}\right), \quad u^{\epsilon}_{x}\in L^{2}\left(0,T; L^{2}\right).
\] 
Moreover, these bounds with the equation 
\[
 \theta^{\epsilon}_{t} =-u^{\epsilon}\theta^{\epsilon}_x-\Lambda^{\gamma}\theta^{\epsilon} +\epsilon \theta^{\epsilon}_{xx} =-\left(u^{\epsilon}\theta^{\epsilon}\right)_x+u^{\epsilon}_{x}\theta^{\epsilon} -\Lambda^{\gamma}\theta^{\epsilon} +\epsilon \theta^{\epsilon}_{xx}
\]
we also have that  
\[
\theta^{\epsilon}_{t} \in L^{\frac{4}{3}}\left(0,T; H^{-2}\right).
\]

We now extract a subsequence of $\left(\theta^{\epsilon}\right)$, using the same index $\epsilon$ for simplicity, and a function $\theta \in \mathcal{C}_T$ and $u=(1-\partial_{xx})^{-\alpha}\theta$ such that 
\begin{equation}\label{convergence of approximate solutions 2}
\begin{split}
& \theta^{\epsilon} \stackrel{\star}{\rightharpoonup} \theta \quad \text{in} \quad L^{\infty}\left(0,T; L^{p}\right) \quad \text{for all $p\in (1,\infty)$},\\
& \theta^{\epsilon} \rightharpoonup \theta \quad \text{in} \quad L^{2}\left(0,T; H^{\frac{\gamma}{2}}\right),\\
& \theta^{\epsilon}\rightarrow \theta \quad \text{in} \quad L^{2}\left(0,T; L^{p}_{\text{loc}}\right) \quad \text{for all $p\in (1,\infty)$},\\
& u^{\epsilon}_{x} \rightharpoonup u_{x} \quad \text{in} \quad L^{2}\left(0,T; L^{2}\right),
\end{split}
\end{equation}
where we use Lemma \ref{lemma:2.1} to obtain the strong convergence. 

We now multiply (\ref{regularized equation 2}) by a test function $\psi\in \mathcal{C}^{\infty}_{c}\left([0,T)\times\mathbb{R}\right)$ and integrate over $\mathbb{R}$. Then,
\[
\int^{T}_{0}\int  \left[\theta^{\epsilon} \psi_{t} +u^{\epsilon}\theta^{\epsilon}\psi_{x}+u^{\epsilon}_{x}\theta^{\epsilon} \psi -\theta^{\epsilon}\Lambda^{\gamma}\psi+\epsilon \theta^{\epsilon}\psi_{xx}\right]  dxdt =\int \theta^{\epsilon}_{0}(x)\psi(0,x)dx.
\]
By Lemma \ref{lemma:2.1} with
\[
X_0=L^2\left(0,T; H^{\frac{\gamma}{2}}\right),\quad X=L^2\left(0,T; L^2_{\text{loc}}\right),\quad X_1=L^2\left(0,T; H^{-2}\right)
\]
and using (\ref{convergence of approximate solutions 2}), we can pass to the limit to obtain that 
\begin{equation} 
\begin{split}
\int^{T}_{0}\int \left[\theta \psi_{t} +u\theta \psi_{x}+u_{x}\theta \psi -\theta \Lambda^{\gamma}\psi\right]  dxdt =\int \theta_{0}(x)\psi(0,x)dx
\end{split}
\end{equation}
This completes the proof.

\textcolor{black}{ 
%%%%%%%%%%%%%%%
\subsection{Uniqueness}
%%%%%%%%%%%%%%%
To show the uniqueness of a weak solution, we consider the equation of $\theta=\theta_{1}-\theta_{2}$ given by 
\eqn \label{eq:6.77}
\theta_{t}+\Lambda^{\gamma} \theta=-u_{1}\theta_{x}+ u\theta_{2x},  \quad \theta(0,x)=0.  
\een
We multiply (\ref{eq:6.77}) by $\theta$ to obtain 
\begin{equation*}
\begin{split}
\frac{1}{2}\left\|\theta\right\|^{2}_{L^{2}}+ \left\|\Lambda^{\frac{\gamma}{2}}\theta \right\|^{2}_{L^{2}}&=-\int u_{1}\theta_{x}\theta dx +\int u\theta_{2x}\theta dx =\frac{1}{2}\int u_{1x}\theta^{2}dx+\int u\theta_{2x}\theta dx\\
& \leq C\|u_{1x}\|_{L^{2}}\|\theta\|_{L^{2}}\|\theta\|_{L^{\infty}} +C\|\theta_{2}\|_{\dot{H}^{1-\frac{\gamma}{2}}} \|u\|_{L^{\infty}} \|\theta\|_{\dot{H}^{\frac{\gamma}{2}}} +C\|\theta_{2}\|_{\dot{H}^{1-\frac{\gamma}{2}}} \|u\|_{\dot{H}^{\frac{\gamma}{2}}}\|\theta\|_{L^{\infty}} \\
& \leq C\|u_{1x}\|_{L^{2}}\|\theta\|_{L^{2}}\|\theta\|_{H^{\frac{\gamma}{2}}} +C\|\theta_{2}\|_{\dot{H}^{1-\frac{\gamma}{2}}} \|u\|_{L^{\infty}} \|\theta\|_{\dot{H}^{\frac{\gamma}{2}}} +C\|\theta_{2}\|_{\dot{H}^{1-\frac{\gamma}{2}}}\|u\|_{\dot{H}^{\frac{\gamma}{2}}} \|\theta\|_{H^{\frac{\gamma}{2}}} \\
& \leq C\left(\|u_{1x}\|_{L^{2}}+ \|u_{1x}\|^{2}_{L^{2}}\right)\|\theta\|^{2}_{L^{2}} +C\|\theta_{2}\|^{2}_{\dot{H}^{1-\frac{\gamma}{2}}} \|u\|^{2}_{L^{\infty}}\\
& +C\left(\|\theta_{2}\|_{\dot{H}^{1-\frac{\gamma}{2}}}+ \|\theta_{2}\|^{2}_{\dot{H}^{1-\frac{\gamma}{2}}} \right)\|u\|^{2}_{\dot{H}^{\frac{\gamma}{2}}}+\frac{1}{2}\left\|\Lambda^{\frac{\gamma}{2}}\theta \right\|^{2}_{L^{2}},
\end{split}
\end{equation*}
where we use $\gamma>1$ to control $\|\theta\|_{L^{\infty}}$. Since $\gamma>1$ and $\alpha>\frac{1}{4}\ge \frac{1}{2}-\frac{\gamma}{4}$, we also have
\[
\|u_{1x}\|_{L^{2}}\leq C \|\theta_{1}\|_{H^{\frac{\gamma}{2}}}, \quad \|u\|_{L^{\infty}} \leq C \|\theta\|_{L^{2}}, \quad \|\theta_{2}\|_{\dot{H}^{1-\frac{\gamma}{2}}} \leq C\|\theta_{2}\|_{H^{\frac{\gamma}{2}}}.
\]
Hence
\[
\left\|\theta\right\|^{2}_{L^{2}}+ \left\|\Lambda^{\frac{\gamma}{2}}\theta \right\|^{2}_{L^{2}}\leq C \left(1+\|\theta_{1}\|^{2}_{H^{\frac{\gamma}{2}}} + \|\theta_{2}\|^{2}_{H^{\frac{\gamma}{2}}} \right)\left\|\theta\right\|^{2}_{L^{2}}.
\]
Since the quantities in the parenthesis are integrable in time and $\theta(0,x)=0$,  we conclude that $\theta=0$ in $L^{2}$ and thus a weak solution is unique. This completes the proof of Theorem \ref{GW 2}.}

%%%%%%%%%%%%%%%%%%%%%%%%%%%%%
\section{Proof of Theorem \ref{GW infinite energy 2}}
%%%%%%%%%%%%%%%%%%%%%%%%%%%%%

%%%%%%%%%%%%%%%%
\subsection{A priori estimate}
%%%%%%%%%%%%%%%%
We consider the equation
\eqn \label{infinite equation 2}
\theta_{t}+u\theta_{x}+\Lambda\theta=0, \quad u=(1-\partial_{xx} )^{-\frac{1}{4}} \theta
\een
We multiply (\ref{infinite equation 2}) by $\theta w$ and integrate in $x$. Then, 
\begin{equation*}
\begin{split}
\frac{1}{2}\frac{d}{dt} \left\|\theta\right\|^{2}_{L^{2}(wdx)}+ \left\|\Lambda^{\frac{1}{2}}\theta\right\|^{2}_{L^{2}(wdx)}&=-\int  u\theta_{x} \theta wdx -\int \Lambda^{\frac{1}{2}}\theta \left[\Lambda^{\frac{1}{2}},w\right]\theta dx\\
& =\frac{1}{2}\int u_{x}\theta^{2} wdx +\frac{1}{2}\int u\theta^{2} w_{x} dx -\int \Lambda^{\frac{1}{2}}\theta \left[\Lambda^{\frac{1}{2}},w\right]\theta dx.
\end{split}
\end{equation*}
As in the proof of Theorem \ref{GW infinite energy}, we bound the commutator term as 
\[
\int \Lambda^{\frac{1}{2}}\theta \left[\Lambda^{\frac{1}{2}},w\right]\theta dx\leq \frac{1}{4}  \left\|\Lambda^{\frac{1}{2}}\theta\right\|^{2}_{L^{2}(wdx)} +  C  \left\|\theta\right\|^{2}_{L^{2}(wdx)}.
\]
To estimate terms involving $u$, we use $\theta=(1-\partial_{xx})^{\frac{1}{4}}u$ and (\ref{weighted Sobolev space}) to obtain that 
\begin{equation*}
\begin{split}
\int u_{x}\theta^{2} wdx + \int u\theta^{2} w_{x}dx &\leq C \|\theta_{0}\|_{L^{\infty}}\left(\left\|u\right\|_{L^{2}(wdx)}+ \left\|u_{x}\right\|_{L^{2}(wdx)}\right) \left\|\theta\right\|_{L^{2}(wdx)}\\
& \leq C \|\theta_{0}\|_{L^{\infty}}\left(\left\|\theta\right\|_{L^{2}(wdx)} + \left\|\Lambda^{\frac{1}{2}}\theta\right\|_{L^{2}(wdx)} \right) \left\|\theta\right\|_{L^{2}(wdx)}\\
&\leq C \left( \|\theta_{0}\|_{L^{\infty}}+  \|\theta_{0}\|^{2}_{L^{\infty}}\right)  \left\|\theta\right\|^{2}_{L^{2}(wdx)} +\frac{1}{4}  \left\|\Lambda^{\frac{1}{2}}\theta\right\|^{2}_{L^{2}(wdx)}.
\end{split}
\end{equation*}

Collecting all terms together, we obtain that 
\[
\frac{d}{dt}\left\|\theta\right\|^{2}_{L^{2}(wdx)}+ \left\|\Lambda^{\frac{1}{2}}\theta\right\|^{2}_{L^{2}(wdx)} \leq C\left(\|\theta_{0}\|_{L^{\infty}}+ \|\theta_{0}\|^{2}_{L^{\infty}}\right)\left\|\theta\right\|^{2}_{L^{2}(wdx)}
\]
and hence that 
\eqn\label{H1 weighted bound 1}
\left\|\theta(t)\right\|^{2}_{L^{2}(wdx)}+\int^{t}_{0} \left\|\Lambda^{\frac{1}{2}}\theta(s)\right\|^{2}_{L^{2}(wdx)} \ ds \leq \left\|\theta_{0}\right\|^{2}_{L^{2}(wdx)} \exp\left[C\left(\|\theta_{0}\|_{L^{\infty}}+ \|\theta_{0}\|^{2}_{L^{\infty}}\right)t\right].
\een

We next multiply (\ref{infinite equation 2}) by $-(\theta_{x} w)_{x}$ and integrate in $x$. Then, 
\begin{equation*}
\begin{split}
&\frac{1}{2}\frac{d}{dt} \left\|\theta_{x}\right\|^{2}_{L^{2}(wdx)}+ \left\|\Lambda^{\frac{3}{2}}\theta\right\|^{2}_{L^{2}(wdx)}=\int  u\theta_{x} (\theta_{x} w)_{x} dx +\int \Lambda^{\frac{3}{2}}\theta \left[\Lambda^{\frac{1}{2}},w\right]\Lambda \theta dx\\
& =-\frac{1}{2}\int u_{x}(\theta_{x})^{2} wdx +\frac{1}{2}\int u(\theta_{x})^{2} w_{x} dx -\int \Lambda^{\frac{1}{2}}\theta \left[\Lambda^{\frac{1}{2}},w\right]\theta dx +\int \Lambda \theta \theta_{x}w_{x}dx.
\end{split}
\end{equation*}
Following the computation in \cite{Lazar}, 
\begin{equation*}
\begin{split}
\frac{d}{dt} \left\|\theta_{x}\right\|^{2}_{L^{2}(wdx)}+ \left\|\Lambda^{\frac{3}{2}}\theta\right\|^{2}_{L^{2}(wdx)}&\leq C \|u\|_{H^{1}(wdx)} \left\|\theta_{x}\right\|^{2}_{L^{4}(wdx)} +C  \left\|\theta_{x}\right\|^{2}_{L^{2}(wdx)} +\frac{1}{4}  \left\|\Lambda^{\frac{3}{2}}\theta\right\|^{2}_{L^{2}(wdx)} \\
& \leq C \|\theta\|_{H^{\frac{1}{2}}(wdx)} \left\|\theta_{x}\right\|^{2}_{L^{4}(wdx)} +C  \left\|\theta_{x}\right\|^{2}_{L^{2}(wdx)} +\frac{1}{4}  \left\|\Lambda^{\frac{3}{2}}\theta\right\|^{2}_{L^{2}(wdx)},
\end{split}
\end{equation*}
where we use the relation in (\ref{weighted Sobolev space}) to bound $u$ in terms of $\theta$. Since
\begin{equation}
\begin{split}
\|\theta\|_{H^{\frac{1}{2}}(wdx)}  \left\|\theta_{x}\right\|^{2}_{L^{4}(wdx)} &\leq C\|\theta\|_{H^{\frac{1}{2}}(wdx)} \|\theta_{x}\|_{L^{2}(wdx)}  \left\|\Lambda^{\frac{3}{2}}\theta\right\|_{L^{2}(wdx)} \\
&\leq C\|\theta\|^{2}_{H^{\frac{1}{2}}(wdx)} \|\theta_{x}\|^{2}_{L^{2}(wdx)}+\frac{1}{4}\left\|\Lambda^{\frac{3}{2}}\theta\right\|^{2}_{L^{2}(wdx)}
\end{split}
\end{equation}
by (\ref{GN weight}), we obtain that 
\eqn \label{H1 weighted bound 2}
\frac{d}{dt} \left\|\theta_{x}\right\|^{2}_{L^{2}(wdx)}+ \left\|\Lambda^{\frac{3}{2}}\theta\right\|^{2}_{L^{2}(wdx)}\leq C\left(1+ \|\theta\|^{2}_{H^{\frac{1}{2}}(wdx)} \right)  \left\|\theta_{x}\right\|^{2}_{L^{2}(wdx)}.
\een
Integrating in time (\ref{H1 weighted bound 2}) and using (\ref{H1 weighted bound 1}), we obtain that
\begin{equation}\label{H1 weighted bound 3}
\begin{split}
&\left\|\theta_{x}(t)\right\|^{2}_{L^{2}(wdx)}+ \int^{t}_{0}\left\|\Lambda^{\frac{3}{2}}\theta(s)\right\|^{2}_{L^{2}(wdx)} \\
&\leq \left\|\theta_{0x}\right\|^{2}_{L^{2}(wdx)}\exp\left[\int^{t}_{0}C\left(1+ \|\theta(s)\|^{2}_{H^{\frac{1}{2}}(wdx)} \right) \ ds \right]\\
& \leq \left\|\theta_{0x}\right\|^{2}_{L^{2}(wdx)}\exp\left[Ct + \left\|\theta_{0}\right\|^{2}_{L^{2}(wdx)} \exp\left[C\left(\|\theta_{0}\|_{L^{\infty}}+ \|\theta_{0}\|^{2}_{L^{\infty}}\right)t\right]\right].
\end{split}
\end{equation}
By (\ref{H1 weighted bound 1}) and (\ref{H1 weighted bound 3}), we finally obtain that 
\begin{equation}\label{H1 weighted bound 4}
\begin{split}
&\left\|\theta(t)\right\|^{2}_{H^{1}(wdx)}+ \int^{t}_{0}\left\|\Lambda^{\frac{1}{2}}\theta(s)\right\|^{2}_{H^{1}(wdx)} \ ds \\
& \leq C\left\|\theta_{0}\right\|^{2}_{H^{1}(wdx)}\exp\left[Ct + \left\|\theta_{0}\right\|^{2}_{L^{2}(wdx)} \exp\left[C\left(\|\theta_{0}\|_{L^{\infty}}+ \|\theta_{0}\|^{2}_{L^{\infty}}\right)t\right]\right].
\end{split}
\end{equation}

%%%%%%%%%%%%%%%%%%%%%%%%%%
\subsection{Approximation and passing to limit}
%%%%%%%%%%%%%%%%%%%%%%%%%%
Since $\theta$ is more regular than a solution in Theorem \ref{GW infinite energy}, we can follow the procedure in the proof of Theorem \ref{GW infinite energy}.

%%%%%%%%%%%%%%
\subsection{Uniqueness}
%%%%%%%%%%%%%%
To show the uniqueness of a weak solution, we consider the equation of $\theta=\theta_{1}-\theta_{2}$ given by 
\eqn \label{eq:2.777}
\theta_{t}+\Lambda \theta=-u_{1}\theta_{x}+ u\theta_{2x},  \quad \theta(0,x)=0.  
\een
We multiply $w\theta$ to (\ref{eq:2.777}) and integrate over $\mathbb{R}$. Then, 
\begin{equation*} 
 \begin{split}
& \frac{1}{2}\frac{d}{dt}\left\|\theta\right\|^{2}_{L^{2}(wdx)}+\left\|\Lambda^{\frac{1}{2}}\theta \right\|^{2}_{L^{2}(wdx)}  = \int \left[-u_{1}\theta_{x}+ u\theta_{2x}\right]\theta wdx -\int \Lambda^{\frac{1}{2}}\theta \left[\Lambda^{\frac{1}{2}},w\right]\theta dx\\
 & =\frac{1}{2}\int u_{1x}\theta^{2}w dx +\frac{1}{2}\int u_{1}\theta^{2}w_{x}dx +\int u\theta_{2x}\theta wdx  -\int \Lambda^{\frac{1}{2}}\theta \left[\Lambda^{\frac{1}{2}},w\right]\theta dx .
 \end{split}
\end{equation*}
As before, the last term is bounded by
\[
\int \Lambda^{\frac{1}{2}}\theta \left[\Lambda^{\frac{1}{2}},w\right]\theta dx \leq \frac{1}{2} \left\|\Lambda^{\frac{1}{2}}\theta\right\|^{2}_{L^{2}(wdx)} +  C \left\|\theta\right\|^{2}_{L^{2}(wdx)}.
\]
The first three terms in the right-hand side are easily bounded by
\[
C\left(\left\|\theta_{2x}\right\|_{L^{2}(wdx)}+\left\|\theta_{1}\right\|_{H^{\frac{1}{2}}(wdx)} \right) \left(\|\theta\|^{2}_{L^{4}(wdx)} + \|\theta\|_{L^{4}(wdx)} \|u\|_{L^{4}(wdx)}\right).
\]
By (\ref{GN weight}), we obtain that 
\begin{equation*} 
\begin{split}
\frac{d}{dt}\left\|\theta\right\|^{2}_{L^{2}(wdx)}+\left\|\Lambda^{\frac{1}{2}}\theta \right\|^{2}_{L^{2}(wdx)}  \leq C \left(1+ \left\|\theta_{2}\right\|^{2}_{H^{1}(wdx)}+ \left\|\theta_{1}\right\|^{2}_{H^{\frac{1}{2}}(wdx)} \right) \|\theta\|^{2}_{L^{2}(wdx)} +\frac{1}{2}\left\|\Lambda^{\frac{1}{2}}\theta \right\|^{2}_{L^{2}(wdx)}.
\end{split}
\end{equation*}
Since 
\[
\theta_{2} \in L^{2}\left(0,T: H^{1}(wdx)\right),\quad \theta_{1} \in L^{2}\left(0,T: H^{\frac{1}{2}}(wdx)\right),
\]
we conclude that $\theta=0$ in $L^{2}(wdx)$ and thus a weak solution is unique. This completes the proof of Theorem \ref{GW infinite energy 2}.

%%%%%%%%%%%%%%%%%%%%%%%%%
\section{Proof of Theorem \ref{GW Wiener2}}
%%%%%%%%%%%%%%%%%%%%%%%%%

Taking the Fourier transform of \eqref{transport equation 2}, we have that
\begin{align*}
\pat|\hat{\theta}(\xi)|&=-\text{Re}\left[\int_{\RR}\frac{1}{(1+|\zeta|^2)^{\alpha}}\hat{\theta}(\zeta)i(\xi-\zeta))\hat{\theta}(\xi-\zeta)d\zeta\frac{\bar{\hat{\theta}}(\xi)}{|\hat{\theta}(\xi)|}\right]\frac{1}{\sqrt{2\pi}}-|\xi||\hat{\theta}|.
\end{align*}
Consequently, ignoring the factor $\frac{1}{(1+|\zeta|^2)^{\alpha}}$, we follow the proof of Theorem \ref{GW Wiener} with $\delta=0$ and the smallness condition (\ref{initial data in Wiener 2}) to obtain that 
\begin{equation}
\|\theta(t)\|_{A^1}+\left(1-\frac{\sqrt{2}\|\theta_{0}\|_{A^0}}{\sqrt{\pi}}\right)\int_0^t \|\theta_{x}(s)\|_{A^1}ds\leq \|\theta_0\|_{A^1}
\end{equation}
for all $t\ge 0$. We also follow the proof of Theorem \ref{GW Wiener}  to obtain a unique weak solution via the approximation procedure. This completes the proof.

%%%%%%%%%%%%%%%
\section*{Acknowledgments}
%%%%%%%%%%%%%%%
H.B. was supported by NRF-2015R1D1A1A01058892. R.G.B. is funded by the LABEX MILYON (ANR-10-LABX-0070) of Universit\'e de Lyon, within the program ``Investissements d'Avenir'' (ANR-11-IDEX-0007) operated by the French National Research Agency (ANR). Both O. L. and R.G.B. were partially supported by the Grant MTM2014-59488-P from the former Ministerio de Econom\'ia y Competitividad (MINECO, Spain). O.L. was partially supported by the Marie-Curie Grant, acronym: TRANSIC, from the FP7-IEF program.

%%%%%%%%%%%%%%%
\bibliographystyle{abbrv}
%\bibliography{bibliografia}

\end{document}